\documentclass[11pt,a4paper]{article}

\usepackage[T1]{fontenc}
\usepackage[utf8]{inputenc}
\usepackage{lmodern}
\usepackage[margin=1in]{geometry}
\usepackage{amsmath,amssymb,amsthm,mathtools}
\usepackage{graphicx}
\usepackage{microtype}
\usepackage[hidelinks]{hyperref}
\usepackage{newtxtext}
\usepackage{newtxmath}
\usepackage{indentfirst}
\usepackage{tikz-cd}
\usepackage[misc]{ifsym}
\usepackage[title]{appendix}
\theoremstyle{plain}
\newtheorem{theorem}{Theorem}[section]
\newtheorem{lemma}[theorem]{Lemma}
\newtheorem{proposition}[theorem]{Proposition}
\newtheorem{corollary}[theorem]{Corollary}
\theoremstyle{definition}
\newtheorem{definition}[theorem]{Definition}

\newtheorem{remark}[theorem]{Remark}
\numberwithin{equation}{section}

\newcommand{\R}{\mathbb{R}}

\newcommand{\HH}{\mathbf H}
\newcommand{\VV}{\mathbf V}
\newcommand{\ff}{\mathbf f}

\newcommand{\dm}{\,d\mu}

\DeclareMathOperator{\Ric}{Ric}

\title{Finite-Time Singularities of the \(L_\beta\)-Gradient Flow in a K\"ahler~Surface}
\author{Weihan Wang}

\makeatletter

\date{} 
\makeatletter
\renewenvironment{abstract}
{%
  \par
  \small
  \noindent\textbf{\abstractname.}\enspace\ignorespaces
}
{%
  \par
  \medskip
}
\makeatother

\begin{document}
\maketitle

\begin{abstract}
We study finite-time singularities of the $L_\beta$-gradient flow of closed
symplectic surfaces in a compact K\"ahler surface, with $\beta>0$ fixed.
We derive the evolution equation for the second fundamental form and obtain
a lower bound on the blow-up rate of its maximum squared norm.
Combining this bound with a local curvature propagation estimate and the
weighted monotonicity formula, we rule out Type~I singularities.
\end{abstract}



\section{Introduction}\label{sec:introduction}

Let $(M,\omega,J)$ be a K\"ahler surface with metric $\bar g$.
For a smooth immersion $F:\Sigma\to M$ of a closed oriented surface,
its K\"ahler angle $\alpha:\Sigma\to[0,\pi]$ is defined by
\[
F^*\omega=\cos\alpha\,d\mu_\Sigma,
\]
where $d\mu_\Sigma$ is the area form of $F^*\bar g$
(see \cite{CW}).
The immersion is called symplectic if $\cos\alpha>0$ everywhere,
Lagrangian if $\cos\alpha\equiv0$, and holomorphic if $\cos\alpha\equiv1$.

By Wirtinger's inequality, closed holomorphic curves in a K\"ahler surface
minimize area in their homology classes and are therefore stable symplectic
minimal surfaces~\cite{HL}.
Wolfson~\cite{W} showed that every closed symplectic minimal surface in a
K\"ahler--Einstein surface with positive scalar curvature is holomorphic.
Thus, in this setting, one may seek holomorphic curves by constructing
closed symplectic minimal surfaces.

A natural way to seek such minimal surfaces is to evolve a closed
symplectic surface by mean curvature flow, the negative gradient flow
of the area functional.
Chen and Li~\cite{CL} and Wang~\cite{Wang} showed that in a
K\"ahler--Einstein surface the evolving surface remains symplectic
for as long as the flow is smooth.
This evolution is called symplectic mean curvature flow.

In a K\"ahler--Einstein surface of negative scalar curvature, however,
a symplectic minimal surface can be nonholomorphic.
Arezzo~\cite{Arezzo} constructed strictly stable embedded two-spheres
with this property.
Since each of these spheres has vanishing mean curvature, the mean
curvature flow starting from it is stationary.

To incorporate the K\"ahler angle into the variational problem, Han, Li, and
Sun~\cite{HLSD1} introduced the family of functionals
\[
L_\beta(\Sigma)=\int_\Sigma \frac{1}{\cos^ \beta\alpha}\,d\mu_\Sigma,
\qquad \beta\geq0,
\]
on symplectic surfaces in a K\"ahler surface.
At $\beta=0$, the functional reduces to area, and its critical points
in this class are precisely the symplectic minimal surfaces.
The authors proposed using the continuity method to follow a family
of critical surfaces from $\beta=0$ in search of a holomorphic limit
as $\beta\to\infty$.
For fixed $\beta$, the critical points of $L_\beta$ are called
$\beta$-symplectic critical surfaces. Their Euler--Lagrange equation is
\[
\cos^3\alpha\,\HH
-\beta\bigl(J(J\nabla\cos\alpha)^\top\bigr)^\perp=0,
\]
where $\HH$ is the mean curvature vector, $\nabla$ is the gradient
with respect to the induced metric, and $(\cdot)^\top$ and
$(\cdot)^\perp$ denote the tangential and normal components,
respectively. Every holomorphic curve satisfies this equation for
all $\beta\geq0$.

To seek critical points of $L_\beta$ for fixed $\beta\geq0$, Han, Li,
and Sun~\cite{HLSflow} considered its negative gradient flow
\begin{equation}\label{eq:gradient-flow}
\begin{cases}
\displaystyle
\frac{\partial F}{\partial t}
=
\frac{\cos^3\alpha\,\HH
      -\beta\bigl(J(J\nabla\cos\alpha)^\top\bigr)^\perp}
     {\cos\alpha\bigl(\cos^2\alpha+\beta\sin^2\alpha\bigr)},
\\[2mm]
F(\cdot,0)=F_0,
\end{cases}
\end{equation}
where $F_0:\Sigma\to M$ is a smooth symplectic immersion.
We denote the immersed surface at time $t$ by $\Sigma_t$.
The stationary solutions are precisely the $\beta$-symplectic critical
surfaces. When $\beta=0$, the evolution equation reduces to the mean
curvature flow equation $\partial_t F=\HH$.
For $\beta>0$, the $L_\beta$-gradient flow preserves symplecticity
in any compact K\"ahler surface for as long as the solution remains
smooth~\cite[Corollary~2.4]{HLSflow}.

If the flow \eqref{eq:gradient-flow} develops a singularity in finite
time, a natural question is how fast the second fundamental form of
$\Sigma_t$ blows up.
For $\beta=0$, Chen and Li~\cite{CL} and Wang~\cite{Wang} ruled out Type~I singularities for symplectic mean curvature flow
in compact K\"ahler--Einstein surfaces.
For $\beta>0$, Han, Li, and Sun~\cite{HLSflow} established a weighted
monotonicity formula and an $\varepsilon$-regularity theorem for the
$L_\beta$-gradient flow. They asked whether Type~I
singularities can also be ruled out in this case.
We answer this question affirmatively.

\begin{theorem}\label{thm:main}
Fix $\beta>0$, and let $F_0:\Sigma\to M$ be a smooth symplectic
immersion of a closed oriented surface into a compact K\"ahler
surface $M$. Then the $L_\beta$-gradient flow
\eqref{eq:gradient-flow} starting from $F_0$ develops no
finite-time Type~I singularities.
\end{theorem}

We derive the evolution equations for the induced metric and the second fundamental form in Section~\ref{sec:evolution-formulae}
and establish a lower bound on the curvature blow-up rate in
Section~\ref{sec:proof}.
In Section~\ref{sec:special-type-I}, we use this bound to rule out Type~I singularities (Theorem~\ref{thm:main}).
Appendix~\ref{app:submanifold-geometry} collects the submanifold
identities used in the proofs.

\section{Evolution equations for the \texorpdfstring{$L_\beta$}{L\_beta}-gradient flow}
\label{sec:evolution-formulae}

We derive the evolution equations for the induced metric and the
second fundamental form along \eqref{eq:gradient-flow} using the
spacetime bundle formalism of Andrews and Baker~\cite{AB}.
Our curvature conventions and the submanifold identities used below
are collected in Appendix~\ref{app:submanifold-geometry}.

Let $I\subset\mathbb R$ be an interval and $F:\Sigma\times I\to M$
a smooth family of immersions, with $F_t=F(\cdot,t)$.
Equip $F^*TM$ with the pullback metric $\bar g_F$ and the pullback
connection ${}^F\bar\nabla$.
Let $t:\Sigma\times I\to I$ denote the projection onto the time factor.
Via $F_*=dF$, we identify the spatial tangent bundle
$\mathcal H=\ker dt$ with the subbundle of $F^*TM$ whose fiber at
$(x,t)$ is the tangent plane $d(F_t)_x(T_x\Sigma)$.
Its orthogonal complement $\mathcal N$ is the normal bundle of the
time slices, giving the orthogonal splitting
\[
F^*TM=\mathcal H\oplus\mathcal N.
\]
We denote the orthogonal projections onto $\mathcal H$ and
$\mathcal N$ by $\pi$ and $\pi^\perp$, respectively, and the
inclusion $\mathcal N\hookrightarrow F^*TM$ by $\iota$.

The restrictions of $\bar g_F$ to $\mathcal H$ and $\mathcal N$
give the bundle metrics $g$ and $g^\perp$, respectively.
For a vector field $X$ on $\Sigma\times I$ and sections
$U$ of $\mathcal H$ and $\xi$ of $\mathcal N$, the induced
connections are given by
\[
\nabla_X U=\pi\bigl({}^F\bar\nabla_X(F_*U)\bigr),
\qquad
\nabla_X^\perp\xi=\pi^\perp\bigl({}^F\bar\nabla_X(\iota\xi)\bigr).
\]
On each time slice, their spatial restrictions are the Levi-Civita
connection of $g(t)=F_t^*\bar g$ and the normal connection of $F_t$,
respectively.
Both connections preserve their respective metrics in all spacetime
directions. In particular, $\nabla_{\partial_t}g=0$, although the
components $g_{ij}(t)$ in fixed spatial coordinates may vary with $t$.

\begin{definition}\label{def:spacetime-second-fundamental-form}
The spacetime second fundamental form of $\mathcal H$ is the
$\mathcal N$-valued tensor $h^{\mathcal H}$ defined by
\[
h^{\mathcal H}(X,V)
:=\pi^\perp\bigl({}^F\bar\nabla_X(F_*V)\bigr),
\qquad
X\in\Gamma(T(\Sigma\times I)),\quad V\in\Gamma(\mathcal H).
\]
For $U,V\in\Gamma(\mathcal H)$, write $h(U,V):=h^{\mathcal H}(U,V)$.
Then $h$ is the second fundamental form of $F_t$.
\end{definition}

Suppose now that $F$ is a family of symplectic immersions evolving by \eqref{eq:gradient-flow}, and let
$\alpha(\cdot,t)$ denote the K\"ahler angle of $F_t$.
At each $(x,t)\in\Sigma\times I$, choose an orthonormal frame
$(e_1,e_2,\nu_3,\nu_4)$ of $T_{F_t(x)}M$ such that
$(e_1,e_2)$ is a positively oriented basis of $\mathcal H_{(x,t)}$
and the complex structure $J$ satisfies
\begin{equation}\label{eq:J-adapted-frame}
\begin{pmatrix}
Je_1\\ Je_2\\ J\nu_3\\ J\nu_4
\end{pmatrix}
=
\begin{pmatrix}
0           & \cos\alpha  & \sin\alpha  & 0\\
-\cos\alpha & 0           & 0           & -\sin\alpha\\
-\sin\alpha & 0           & 0           & \cos\alpha\\
0           & \sin\alpha  & -\cos\alpha & 0
\end{pmatrix}
\begin{pmatrix}
e_1\\ e_2\\ \nu_3\\ \nu_4
\end{pmatrix}.
\end{equation}
We use $i,j,k,l,p,q\in\{1,2\}$ for tangent indices and
$\gamma,\zeta\in\{3,4\}$ for normal indices.
Repeated indices are summed over unless otherwise stated.

In the adapted frame, write
$h(e_i,e_j)=\sum_{\gamma=3}^4 h_{ij}^\gamma\nu_\gamma$.
At a complex point, $\cos\alpha=\bar g(Je_1,e_2)$ attains its
maximum value of $1$. Hence $\nabla\cos\alpha=0$ there, and
\eqref{eq:gradient-flow} gives $\partial_t F=\HH$.
Away from complex points, set $\alpha_i:=e_i(\alpha)$.
By \cite[Eq.~(2.5)]{HanLi2010},
\[
\alpha_i=-(h_{i1}^{4}+h_{i2}^{3}),\qquad i=1,2.
\]
On this set, define the normal vector field
$\VV:=\alpha_2\nu_3+\alpha_1\nu_4$.

For each $t$, let $E_t:=\{x\in\Sigma:\cos\alpha(x,t)=1\}$ be the
set of complex points and set $\Omega_t:=\Sigma\setminus E_t$.
For global integral identities, extend $\VV$ from $\Omega_t$ by
zero on $E_t$ and denote this extension by $\widetilde{\VV}$.
Then $|\widetilde{\VV}|\leq\sqrt{2}\,|h|$, and the identity
\[
(1-\cos^2\alpha)\widetilde{\VV}
=\frac{\bigl(J(J\nabla\cos\alpha)^\top\bigr)^\perp}{\cos\alpha}
\]
shows that $(1-\cos^2\alpha)\widetilde{\VV}$ is smooth across $E_t$.
We write $\VV$ for the measurable extension whenever no derivative
of $\VV$ is involved. Formulas involving derivatives of $\alpha$
or $\VV$ are understood on $\Omega_t$.

By the preceding identity, the flow equation \eqref{eq:gradient-flow}
takes the form
\begin{equation}\label{eq:gradient-flow-f}
\partial_t F=\ff
:=\frac{\cos^2\alpha\,\HH-\beta\sin^2\alpha\,\VV}
        {\cos^2\alpha+\beta\sin^2\alpha},
\end{equation}
where $\HH=\operatorname{tr}_{g(t)}h$ is the mean curvature vector
and $\ff$ is the normal velocity of $F_t$. 
We write $H^\gamma$, $V^\gamma$, and $f^\gamma$ for the components of $\HH$, $\VV$, and $\ff$ in the normal frame $(\nu_3,\nu_4)$, respectively, with $H^\gamma=h_{11}^\gamma+h_{22}^\gamma$.

\begin{lemma}[Evolution of the metric and area form]
\label{lem:metric-evolution}
In oriented local coordinates $(x^1,x^2)$ on $\Sigma$ fixed in time,
write $g_{ij}=g(\partial_i,\partial_j)$ and
$h_{ij}=h(\partial_i,\partial_j)$.
Along the flow \eqref{eq:gradient-flow-f}, the induced metric and
area form satisfy
\begin{align}
\partial_t g_{ij}
&=-2g^\perp(h_{ij},\ff),
\label{eq:metric-evolution}\\
\partial_t d\mu_g
&=-g^\perp(\HH,\ff)\,d\mu_g.
\label{eq:area-variation}
\end{align}
\end{lemma}

\begin{proof}
Write $F_i=F_*\partial_i$. Since $[\partial_t,\partial_i]=0$,
the pullback torsion identity (Proposition~\ref{prop:pullback-torsion})
gives
\[
{}^F\bar\nabla_{\partial_t}F_i
={}^F\bar\nabla_{\partial_i}(\iota\ff).
\]
The Weingarten relation \eqref{eq:weingarten2} yields
\[
\bar g_F\bigl({}^F\bar\nabla_{\partial_i}(\iota\ff),F_j\bigr)
=-g^\perp(h_{ij},\ff).
\]
Differentiating $g_{ij}=\bar g_F(F_i,F_j)$ in time and using
the symmetry of $h$ gives \eqref{eq:metric-evolution}.

Differentiating $d\mu_g=\sqrt{\det(g_{ij})}\,dx^1\wedge dx^2$
in the same coordinates, we obtain
\[
\partial_t d\mu_g
=\tfrac12 g^{ij}\partial_t g_{ij}\,d\mu_g
=-g^\perp(\HH,\ff)\,d\mu_g,
\]
where the last equality follows from \eqref{eq:metric-evolution}
and $\HH=g^{ij}h_{ij}$.

The induced spacetime connection preserves $g$ and hence also
$g^{-1}$ and the associated area form. Thus
$\nabla_{\partial_t}g=\nabla_{\partial_t}g^{-1}=0$ and
$\nabla_{\partial_t}d\mu_g=0$.
\end{proof}

We next compute the covariant time derivative of $h$ under the same normal variation.

\begin{proposition}\label{prop:second-fundamental-form-evolution}
Along the flow \eqref{eq:gradient-flow-f}, the second fundamental
form satisfies
\begin{equation}\label{eq:evol_h_in_N}
(\nabla_{\partial_t}h)(U,V)
=(\nabla^{\perp 2}\ff)(U,V)
+h\bigl(V,\mathcal W(U,\ff)\bigr)
+\pi^\perp\bigl(\bar R(\iota\ff,F_*U)F_*V\bigr)
\end{equation}
for $U,V\in\Gamma(\mathcal H)$. Here $\mathcal W$ is the Weingarten
map and
\[
(\nabla^{\perp 2}\ff)(U,V)
:=\nabla^\perp_U\nabla^\perp_V\ff-\nabla^\perp_{\nabla_UV}\ff.
\]
In orthonormal tangent and normal frames, with $h_{ij}=h(e_i,e_j)$,
this identity reads
\begin{equation}\label{evolution_eq_second_fundamental_form}
(\nabla_{\partial_t}h)_{ij}
=(\nabla^{\perp 2}\ff)_{ij}
+g^\perp(\ff,h_{ip})h_{pj}
+f^\gamma\bar R_{i\gamma j\zeta}\nu_\zeta,
\end{equation}
where $\bar R_{i\gamma j\zeta}
=\bar R(F_*e_i,\iota\nu_\gamma,F_*e_j,\iota\nu_\zeta)$.
\end{proposition}

\begin{proof}
Since $F_*\partial_t=\iota\ff$ is normal, the mixed space--time
Codazzi identity of Proposition~\ref{timelikeCodazzi} gives
\eqref{eq:evol_h_in_N}.
In the chosen frames, the Weingarten relation \eqref{eq:weingarten2}
and the symmetry of $h$ give
\[
h\bigl(e_j,\mathcal W(e_i,\ff)\bigr)
=g^\perp(\ff,h_{ip})h_{pj}.
\]
Writing $\ff=f^\gamma\nu_\gamma$ and using our curvature convention,
we obtain
\[
\pi^\perp\bigl(\bar R(\iota\ff,F_*e_i)F_*e_j\bigr)
=-f^\gamma\bar R_{\gamma i j\zeta}\nu_\zeta
=f^\gamma\bar R_{i\gamma j\zeta}\nu_\zeta.
\]
Substituting these two expressions into \eqref{eq:evol_h_in_N}
yields \eqref{evolution_eq_second_fundamental_form}.
\end{proof}

Taking the inner product of \eqref{eq:evol_h_in_N} with $2h$ gives the evolution equation for $|h|^2$. To express it in terms of $\HH$ and $\VV$, write $\ff=\lambda\HH+\mu\VV$, where
\[
\lambda:=\frac{\cos^2\alpha}{\cos^2\alpha+\beta\sin^2\alpha},
\qquad
\mu:=\lambda-1.
\]
For a normal vector $\xi$ and spatial tangent vectors $U,V$, set
\[
\mathcal A_\xi(U,V):=h\bigl(V,\mathcal W(U,\xi)\bigr),
\qquad
\mathcal R_\xi(U,V):=
\pi^\perp\bigl(\bar R(\iota\xi,F_*U)F_*V\bigr).
\]

\begin{proposition}\label{prop:evolution-norm-h-invariant}
Along the flow \eqref{eq:gradient-flow-f}, the squared norm of the
second fundamental form satisfies
\begin{equation}\label{eq:norm-h-general-variation}
\partial_t|h|^2
=2\bigl\langle
  \nabla^{\perp 2}\ff+\mathcal A_{\ff}+\mathcal R_{\ff},h
\bigr\rangle.
\end{equation}
At noncomplex points, this identity takes the form
\begin{align}
(\partial_t-\lambda\Delta)|h|^2
={}&-2\lambda|\nabla h|^2
    -2\lambda\langle\mathcal S,h\rangle
\notag\\
&+2\mu\bigl\langle
  \nabla^{\perp 2}\VV+\mathcal A_{\VV}+\mathcal R_{\VV},h
\bigr\rangle
\notag\\
&+2(\nabla^2\lambda)_{ij}\,g^\perp(\HH+\VV,h_{ij})
\notag\\
&+4e_i(\lambda)\,
  g^\perp\bigl(\nabla^\perp_{e_j}(\HH+\VV),h_{ij}\bigr).
\label{eq:evolution-norm-h-invariant}
\end{align}
Here $(\nabla^2\lambda)_{ij}=(\nabla^2\lambda)(e_i,e_j)$ in a local
orthonormal tangent frame, $\Delta=\operatorname{tr}_g\nabla^2$,
and $\mathcal S=\mathcal S(h,\bar R,\bar\nabla\bar R)$ is the
remainder in the traced Simons identity \eqref{traced-Simons-identity}.
\end{proposition}

\begin{proof}
Metric compatibility gives
$\partial_t|h|^2=2\langle\nabla_{\partial_t}h,h\rangle$.
Substituting \eqref{eq:evol_h_in_N} yields
\eqref{eq:norm-h-general-variation}.

At a noncomplex point, the product rule, $\mu=\lambda-1$, and
the symmetry of $h$ give
\begin{align*}
\langle\nabla^{\perp 2}\ff,h\rangle
={}&\lambda\langle\nabla^{\perp 2}\HH,h\rangle
   +\mu\langle\nabla^{\perp 2}\VV,h\rangle\\
&+(\nabla^2\lambda)_{ij}\,g^\perp(\HH+\VV,h_{ij})\\
&+2e_i(\lambda)\,
  g^\perp\bigl(\nabla^\perp_{e_j}(\HH+\VV),h_{ij}\bigr).
\end{align*}
Since $\mathcal A_\xi$ and $\mathcal R_\xi$ are linear in $\xi$,
the terms with coefficient $2\lambda$ in
\eqref{eq:norm-h-general-variation} combine to give
\[
2\lambda\bigl\langle
\nabla^{\perp 2}\HH+\mathcal A_{\HH}+\mathcal R_{\HH},h
\bigr\rangle
=2\lambda\langle\Delta h-\mathcal S,h\rangle
\]
by \eqref{traced-Simons-identity}. Finally, the identity
\[
2\langle\Delta h,h\rangle=\Delta|h|^2-2|\nabla h|^2
\]
gives \eqref{eq:evolution-norm-h-invariant}.
\end{proof}

We trace \eqref{eq:evol_h_in_N} and pair the result with $2\HH$ to obtain the evolution equation for $|\HH|^2$.

\begin{proposition}
\label{prop:evolution-norm-H-invariant}
Along the flow \eqref{eq:gradient-flow-f}, the squared norm of
the mean curvature vector satisfies
\begin{equation}\label{eq:norm-H-general-variation}
\partial_t|\HH|^2
=2g^\perp\bigl(
  \Delta^\perp\ff
  +\operatorname{tr}_g(\mathcal A_{\ff}+\mathcal R_{\ff}),\HH
\bigr),
\end{equation}
where $\Delta^\perp=\operatorname{tr}_g\nabla^{\perp 2}$.
At noncomplex points, this identity takes the form
\begin{align}
(\partial_t-\lambda\Delta)|\HH|^2
={}&-2\lambda|\nabla^\perp\HH|^2
    +2\mu\,g^\perp(\Delta^\perp\VV,\HH)
\notag\\
&+2(\Delta\lambda)\,g^\perp(\HH+\VV,\HH)
\notag\\
&+4e_i(\lambda)\,
  g^\perp\bigl(\nabla^\perp_{e_i}(\HH+\VV),\HH\bigr)
\notag\\
&+2g^\perp\bigl(
  \operatorname{tr}_g(\mathcal A_{\ff}+\mathcal R_{\ff}),\HH
\bigr).
\label{eq:evolution-norm-H-invariant}
\end{align}
\end{proposition}

\begin{proof}
Since $\HH=\operatorname{tr}_g h$ and $\nabla_{\partial_t}g=0$,
taking the trace of \eqref{eq:evol_h_in_N} gives
\[
\nabla^\perp_{\partial_t}\HH
=\Delta^\perp\ff
+\operatorname{tr}_g(\mathcal A_{\ff}+\mathcal R_{\ff}).
\]
Taking the inner product of both sides with $2\HH$ gives
\eqref{eq:norm-H-general-variation}.

At a noncomplex point, the product rule and $\mu=\lambda-1$ give
\[
\begin{aligned}
\Delta^\perp\ff
={}&\lambda\Delta^\perp\HH+\mu\Delta^\perp\VV+2e_i(\lambda)\nabla^\perp_{e_i}(\HH+\VV)
  +(\Delta\lambda)(\HH+\VV).
\end{aligned}
\]
Substituting this expansion into \eqref{eq:norm-H-general-variation} and using
\[
2g^\perp(\Delta^\perp\HH,\HH)
=\Delta|\HH|^2-2|\nabla^\perp\HH|^2
\]
yields \eqref{eq:evolution-norm-H-invariant}.
\end{proof}

\section{A lower bound for the curvature blow-up rate}
\label{sec:proof}
We estimate the blow-up rate of $\max_{\Sigma_t}|h|^2$ at a finite
singular time of the flow \eqref{eq:gradient-flow}. Throughout this
section, $\beta>0$ is fixed and the initial map $F_0:\Sigma\to M$
is a smooth symplectic immersion of a closed oriented surface into
a compact K\"ahler surface.

The principal-symbol calculation in \cite[Proposition~3.1]{HLSD1}
shows that the flow is parabolic modulo tangential reparametrizations.
In a suitable parametrization, standard parabolic theory yields
short-time existence. For continuation beyond a finite time, we use the following curvature criterion.

\begin{theorem}[Theorem~3.1 of \cite{HLSflow}]
\label{thm:extension}
Let $F:\Sigma\times[0,T)\to M$ be a smooth solution of
\eqref{eq:gradient-flow}, where $0<T<\infty$. If
\[
\sup_{0\leq t<T}\max_{\Sigma_t}|h|^2<\infty,
\]
then $F$ extends smoothly to $[0,T+\epsilon)$ for some $\epsilon>0$.
\end{theorem}

The curvature estimates below require a uniform positive lower bound for $\cos\alpha$ on finite time intervals.

\begin{theorem}[Corollary~2.4 of \cite{HLSflow}]
\label{bound_kahler_angle}
Assume $|\Ric_M|\leq K_1$. Then along any smooth solution of
\eqref{eq:gradient-flow},
\begin{equation}\label{eq:2.3}
\min_{\Sigma_t}\cos\alpha
\geq e^{-K_1t/\beta}\min_{\Sigma_0}\cos\alpha.
\end{equation}
In particular, $\Sigma_t$ remains symplectic for as long as the
solution is smooth.
\end{theorem}

If the maximal existence time $T$ is finite,
Theorem~\ref{thm:extension} implies that $\max_{\Sigma_t}|h|^2$ is unbounded as $t\uparrow T$. 
We now give a lower bound for its blow-up rate.

\begin{theorem}\label{theorem:blowup_lower_rate}
Let $F:\Sigma\times[0,T)\to M$ be a maximal smooth solution of
\eqref{eq:gradient-flow} with $0<T<\infty$. Assume that the initial
immersion $F_0$ is symplectic. Set $v_0:=\min_{\Sigma_0}\cos\alpha$
and $Q(t):=\max_{\Sigma_t}|h|^2$. Then there exists a constant
$\theta=\theta(M,\bar g,J,\beta,v_0,T)>0$ such that
\begin{equation}\label{eq:blowup_lower_rate}
Q(t)\geq\frac{\theta}{2(T-t)},
\qquad 0<T-t\leq\theta/2.
\end{equation}
\end{theorem}

At a noncomplex point where $|h|^2$ attains a spatial maximum,
the term $2\mu\langle\nabla^{\perp 2}\VV,h\rangle$ in
\eqref{eq:evolution-norm-h-invariant} has no definite sign.
To prove Theorem~\ref{theorem:blowup_lower_rate}, we estimate how
long $Q(t)$ stays below a fixed multiple of $1+Q(t_0)$ after time $t_0$.
Unless otherwise stated, the constants below depend only on
$(M,\bar g,J)$, $\beta$, $v_0$, and $T$.
We begin by computing the principal symbol of the normal
linearization of \eqref{eq:gradient-flow}.

Fix $(x_0,t_0)\in\Sigma\times[0,T)$ and let $\mathcal N_{t_0}$
denote the normal bundle of $F_{t_0}$. For
$U\in\Gamma(\mathcal N_{t_0})$, consider the normal variation
$F_\varepsilon(x):=\exp_{F_{t_0}(x)}(\varepsilon U(x))$.
Let $\ff_\varepsilon$ denote the right-hand side of
\eqref{eq:gradient-flow} evaluated at $F_\varepsilon$.
The normal linearization at $F_{t_0}$ is given by
\[
\mathscr L_{t_0}U
:=\pi^\perp_{t_0}
  \left.\frac{D\ff_\varepsilon}{d\varepsilon}\right|_{\varepsilon=0},
\]
where $D/d\varepsilon$ is the ambient covariant derivative along
the variation and $\pi^\perp_{t_0}$ is the orthogonal projection
onto $\mathcal N_{t_0}$.

Since $\bar\nabla\omega=0$, differentiating
$\cos\alpha=\omega(e_1,e_2)$ shows that $\nabla\cos\alpha$ is
linear in $h$. Thus the normal velocity in \eqref{eq:gradient-flow}
is linear in $h$, with coefficients depending smoothly on $F$
and $dF$. Its normal linearization $\mathscr L_{t_0}$ is a
second-order linear differential operator whose principal part
in local coordinates is
\[
\mathcal A^{ij}(\nabla^{\perp 2}U)_{ij},
\]
where each $\mathcal A^{ij}$ is an endomorphism of
$\mathcal N_{t_0}$.

To compute the principal symbol at $x_0$, choose oriented geodesic
coordinates $\psi:B_r(0)\subset\mathbb R^2\to\Sigma$ for $g(t_0)$
with $\psi(0)=x_0$, and set
$e_i:=d(F_{t_0}\circ\psi)_0(\partial_i)$ for $i=1,2$.
Choose a local orthonormal normal frame $(\nu_3,\nu_4)$ such that
$(e_1,e_2,\nu_3,\nu_4)$ is adapted to $J$ at $x_0$ as in
\eqref{eq:J-adapted-frame}.
Writing $U=w^\gamma\nu_\gamma$ and $w=(w^3,w^4)^{\mathsf T}$,
we represent $\mathscr L_{t_0}$ in these coordinates and this frame
by the operator
\[
L_{t_0}w
=\mathcal A^{ij}(y)\partial_i\partial_j w
 +\mathcal B^i(y)\partial_i w
 +\mathcal C(y)w,
\]
where $\mathcal A^{ij},\mathcal B^i,\mathcal C:
B_r(0)\to\operatorname{End}(\mathbb R^2)$ are smooth.

For an operator $L$ of the above form, we use the positive principal
symbol
\[
P_L(y,\xi):=\mathcal A^{ij}(y)\xi_i\xi_j,
\qquad y\in B_r(0),\quad \xi\in\mathbb R^2.
\]
Set $a_0:=\lambda(x_0,t_0)$. In the adapted frame at $x_0$,
the principal part of \eqref{eq:gradient-flow-f} gives
\[
P_{L_{t_0}}(0,\xi)
=
\begin{pmatrix}
a_0\xi_1^2+\xi_2^2 & (1-a_0)\xi_1\xi_2 \\
(1-a_0)\xi_1\xi_2 & \xi_1^2+a_0\xi_2^2
\end{pmatrix}.
\]
Thus, for $\eta=(\eta_3,\eta_4)^{\mathsf T}\in\mathbb R^2$,
\[
\bigl\langle P_{L_{t_0}}(0,\xi)\eta,\eta\bigr\rangle
=a_0|\xi|^2|\eta|^2
 +(1-a_0)(\xi_2\eta_3+\xi_1\eta_4)^2.
\]

By \eqref{eq:2.3}, $\cos\alpha\geq\delta:=e^{-K_1T/\beta}v_0>0$
on $\Sigma\times[0,T)$. The definition of $\lambda$ then gives
\[
0<\kappa:=\frac{\delta^2}{\delta^2+\beta(1-\delta^2)}
\leq a_0\leq 1.
\]
Using $|\xi_2\eta_3+\xi_1\eta_4|\leq|\xi|\,|\eta|$ in the
quadratic form above, we obtain
\[
\kappa|\xi|^2|\eta|^2
\leq \bigl\langle P_{L_{t_0}}(0,\xi)\eta,\eta\bigr\rangle
\leq |\xi|^2|\eta|^2,
\]
uniformly in $(x_0,t_0)\in\Sigma\times[0,T)$.

Freezing the principal coefficients of $L_{t_0}$ at $y=0$ gives
the constant-coefficient operator
$\mathcal A^{ij}(0)\partial_i\partial_j$ on $\mathbb R^2$.
We seek an $L^1$ bound for its heat kernel that is uniform in
$(x_0,t_0)$.

\begin{definition}\label{def:heat_kernel}
Let $L=\mathcal A^{ij}\partial_i\partial_j$ be a constant-coefficient
operator acting on $\mathbb R^2$-valued functions on $\mathbb R^2$.
Assume that $P_L(\xi)$ is symmetric and positive definite for
$\xi\neq0$. For $s>0$, the matrix heat kernel $G_s:\mathbb R^2\to\operatorname{End}(\mathbb R^2)$ of $L$ is defined by
\[
G_s(y):=\frac{1}{(2\pi)^2}
\int_{\mathbb R^2} e^{\mathrm{i}y\cdot\xi}
e^{-sP_L(\xi)}\,d\xi,
\qquad y\in\mathbb R^2,
\]
where $e^{-sP_L(\xi)}$ denotes the matrix exponential.
\end{definition}

\begin{lemma}\label{lem:frozen-kernel-bound}
Let $L$ and $G_s$ be as in Definition~\ref{def:heat_kernel}.
Suppose that
\[
\kappa|\xi|^2 I\leq P_L(\xi)\leq|\xi|^2 I
\qquad\text{for all }\xi\in\mathbb R^2
\]
for some $\kappa>0$. Then there exists a constant
$C_0=C_0(\kappa)\geq1$ such that
\[
\int_{\mathbb R^2}\|G_s(y)\|_{\mathrm{op}}\,dy\leq C_0
\qquad\text{for all }s>0.
\]
\end{lemma}

\begin{proof}
Write $P=P_L$. Since $P$ is homogeneous of degree two,
$G_s(y)=s^{-1}G_1(y/\sqrt{s})$, so it suffices to consider $s=1$.
Set $M(\xi):=e^{-P(\xi)}$, so that $G_1=\mathcal F^{-1}M$.
Since $(1+|y|^2)^{-1}\in L^2(\mathbb R^2)$, the Cauchy--Schwarz
inequality and Plancherel's theorem give
\begin{equation}\label{eq:fourier-L1-estimate}
\int_{\mathbb R^2}\|G_1(y)\|_{\mathrm{op}}\,dy
\leq C\|(1+|y|^2)G_1\|_{L^2}
\leq C\bigl(\|M\|_{L^2}+\|\Delta_\xi M\|_{L^2}\bigr),
\end{equation}
where the matrix-valued $L^2$ norms use the Frobenius norm.

By polarization, the upper bound on the quadratic symbol gives
$\|\partial_{\xi_j}P(\xi)\|_{\mathrm{op}}\leq C|\xi|$ and
$\|\partial_{\xi_j}^2P\|_{\mathrm{op}}\leq C$ with a universal
constant $C$. The lower bound gives
$\|e^{-rP(\xi)}\|_{\mathrm{op}}\leq e^{-r\kappa|\xi|^2}$ for $r\geq0$.
Differentiating the Duhamel formula
\[
\partial_{\xi_j}M
=-\int_0^1 e^{-(1-r)P}(\partial_{\xi_j}P)e^{-rP}\,dr
\]
once more and using these bounds, we obtain
\[
\|M(\xi)\|_{\mathrm F}+\|\Delta_\xi M(\xi)\|_{\mathrm F}
\leq C(1+|\xi|^2)e^{-\kappa|\xi|^2}.
\]
The $L^2(\mathbb R^2)$ norm of the right-hand side depends only
on $\kappa$. The result follows from \eqref{eq:fourier-L1-estimate} and scaling.
Enlarging the constant if necessary, we may take $C_0\geq1$.
\end{proof}

\begin{corollary}\label{cor:frozen-semigroup-bound}
Let $L$ and $G_s$ be as in Lemma~\ref{lem:frozen-kernel-bound}.
Suppose that $w:\mathbb R^2\times(0,S]\to\mathbb R^2$ is a bounded
classical solution of
\[
\left\{
\begin{aligned}
\partial_s w &= Lw
&&\text{on }\mathbb R^2\times(0,S],\\
w(\cdot,s) &\longrightarrow w_0
&&\text{in }\mathcal D'(\mathbb R^2)
\text{ as }s\downarrow0,
\end{aligned}
\right.
\]
where $w_0\in L^\infty(\mathbb R^2;\mathbb R^2)$. Then
$w(\cdot,s)=G_s*w_0$ and
$\|w(\cdot,s)\|_{L^\infty}\leq C_0\|w_0\|_{L^\infty}$
for every $s\in(0,S]$.
Here $C_0$ is the constant given in
Lemma~\ref{lem:frozen-kernel-bound}.
\end{corollary}

\begin{proof}
Fix $s\in(0,S]$ and
$\varphi\in C_c^\infty(\mathbb R^2;\mathbb R^2)$. Set
\[
\Phi(y,r):=\int_{\mathbb R^2}
G_{s-r}(x-y)^{\mathsf T}\varphi(x)\,dx,
\qquad 0\leq r<s.
\]
The Fourier formula for $G_s$ gives
$\partial_r\Phi=-L^*\Phi$, where $L^*$ is the formal adjoint of $L$,
and $\Phi(\cdot,r)\to\varphi$ in $L^1$ as $r\uparrow s$.
Using a spatial cutoff and the rapid decay of $\Phi$, integration
by parts gives
\[
\int_{\mathbb R^2}\langle w(y,s),\varphi(y)\rangle\,dy
=
\int_{\mathbb R^2}\langle w(y,\tau),\Phi(y,\tau)\rangle\,dy
\qquad (0<\tau<s).
\]

Since $w$ is uniformly bounded, the density of $C_c^\infty$ in
$L^1$ extends its distributional convergence at time zero to
convergence against every $L^1$ test function. Together with
$\Phi(\cdot,\tau)\to\Phi(\cdot,0)$ in $L^1$, this allows us to
let $\tau\downarrow0$ and obtain
\[
\begin{aligned}
\int_{\mathbb R^2}\langle w(y,s),\varphi(y)\rangle\,dy
&=\int_{\mathbb R^2}\langle w_0(y),\Phi(y,0)\rangle\,dy\\
&=\int_{\mathbb R^2}
\langle (G_s*w_0)(x),\varphi(x)\rangle\,dx,
\end{aligned}
\]
where the last equality follows from Fubini's theorem.
Since $\varphi$ is arbitrary and both sides of
$w(\cdot,s)=G_s*w_0$ are continuous, this identity holds pointwise.
Finally, Lemma~\ref{lem:frozen-kernel-bound} gives
\[
\|w(\cdot,s)\|_{L^\infty}
\leq
\left(\int_{\mathbb R^2}\|G_s(y)\|_{\mathrm{op}}\,dy\right)
\|w_0\|_{L^\infty}
\leq C_0\|w_0\|_{L^\infty}.
\]
\end{proof}

The following local estimate will be used to compare tangent planes
at different times.

\begin{lemma}[Chord--tangent estimate]\label{lem:chord-tangent}
Let $F:\Sigma\to(U,\bar g)$ be an immersion, where
$U\subset\mathbb R^4$ and $\bar g$ is uniformly equivalent to the
Euclidean metric. For any unit-speed geodesic
$\gamma:[0,r]\to\Sigma$ of $F^*\bar g$, with $r>0$, we have
\[
\left|
\frac{F(\gamma(r))-F(\gamma(0))}{r}
-dF_{\gamma(0)}\bigl(\gamma'(0)\bigr)
\right|
\leq Cr\left(
\sup_{\gamma([0,r])}|h|
+\|\Gamma(\bar g)\|_{C^0(U)}
\right).
\]
Here the left-hand norm is Euclidean, and
$\|\Gamma(\bar g)\|_{C^0(U)}$ is computed using the Euclidean norm of the Christoffel coefficients.
The constant $C$ depends only on the metric equivalence constants.
\end{lemma}

\begin{proof}
Set $V=(F\circ\gamma)'$. The Gauss formula and the geodesic
equation give
$V'=h(\gamma',\gamma')-\Gamma(\bar g)(V,V)$ in coordinates.
Since $|V|_{\bar g}=1$, metric equivalence yields
\[
\sup_{[0,r]}|V'|
\leq C\left(
\sup_{\gamma([0,r])}|h|
+\|\Gamma(\bar g)\|_{C^0(U)}
\right).
\]
The estimate follows from
\[
\left|
\frac{F(\gamma(r))-F(\gamma(0))}{r}-V(0)
\right|
=\left|\frac1r\int_0^r(r-\tau)V'(\tau)\,d\tau\right|
\leq\frac r2\sup_{[0,r]}|V'|.
\]
\end{proof}

We now use these estimates in a rescaling argument to control
the growth of $Q(t)$.

\begin{proposition}\label{prop:curvature-persistence}
There exist constants $L>1$ and $\theta>0$, depending only on
$(M,\bar g,J)$, $\beta$, $v_0$, and $T$, such that for every
$t_0\in[0,T)$,
\[
Q(t)<LK
\qquad\text{for }t_0\leq t<
\min\left\{T,t_0+\frac{\theta}{K}\right\},
\]
where $K:=1+Q(t_0)$.
\end{proposition}

\begin{proof}
Let $\delta>0$ be the uniform lower bound for $\cos\alpha$ obtained
above, let $C_0$ be the constant in
Lemma~\ref{lem:frozen-kernel-bound}, and set 
\[
L=(4C_0)^2.
\]
Suppose that no uniform $\theta$ exists. Then, with
$(M,\bar g,J)$, $\beta$, $v_0$, and $T$ fixed, for each $m$
we can choose a solution and times $t_m\in[0,T)$ and
\[
\tau_m\in[t_m,T)\cap
\left[t_m,t_m+\frac{1}{mK_m}\right),
\qquad K_m:=1+Q(t_m),
\]
such that $Q(\tau_m)\geq LK_m$. Let $b_m\in(t_m,\tau_m]$ be the
first time at which $Q(b_m)=LK_m$, and choose $p_m\in\Sigma$ such
that $|h(p_m,b_m)|^2=LK_m$. Thus
\begin{equation}\label{eq:prop39m-first-hit}
Q(t)\leq LK_m
\qquad (t_m\leq t\leq b_m).
\end{equation}

\smallskip
\noindent\emph{Step 1. Rescaling.}
Set
\[
d_m:=b_m-t_m,
\qquad
\sigma_m:=\sqrt{K_md_m},
\qquad
s:=\frac{t-t_m}{d_m},
\]
and rescale by
\[
\bar g_m:=d_m^{-1}\bar g,
\qquad
F_m(x,s):=F(x,t_m+d_ms).
\]
Let $h_m$, $\ff_m$, and $\alpha_m$ denote the second fundamental
form, the velocity, and the K\"ahler angle of $F_m$ with respect to
$\bar g_m$, and let $g_m(s):=F_m(\cdot,s)^*\bar g_m$. Since
$0<d_m<1/(mK_m)$,
\begin{equation}\label{eq:prop39m-scales}
\sigma_m^2<\frac1m,
\qquad
\frac{d_m}{\sigma_m}=\frac{\sigma_m}{K_m}\leq\sigma_m.
\end{equation}
Hence $\sigma_m\to0$ and $d_m/\sigma_m\to0$.

Constant metric rescaling leaves the connection and the
tangent--normal splitting unchanged, while
$|h_m|=\sqrt{d_m}|h|$. Moreover, if $(e_1,e_2)$ is an oriented
orthonormal tangent frame whose intrinsic covariant derivatives
vanish at the point under consideration, then
$\bar\nabla\omega=0$ gives
\[
d(\cos\alpha)(Z)
=\omega(h(Z,e_1),e_2)+\omega(e_1,h(Z,e_2)).
\]
Thus the right-hand side of \eqref{eq:gradient-flow} depends
smoothly and linearly on $h$ on the set
$\cos\alpha\geq\delta$, including at $\cos\alpha=1$, and
$|\ff|\leq C(\beta,\delta)|h|$. It follows from
\eqref{eq:prop39m-first-hit} that
\begin{equation}\label{eq:prop39m-rescaled-bounds}
\begin{aligned}
\|h_m(\cdot,0)\|_{L^\infty}\leq\sigma_m,
&\qquad
\|h_m\|_{L^\infty(\Sigma\times[0,1])}
\leq\sqrt L\,\sigma_m,\\
|h_m(p_m,1)|=\sqrt L\,\sigma_m,
&\qquad
\|\ff_m\|_{L^\infty(\Sigma\times[0,1])}\leq C\sigma_m.
\end{aligned}
\end{equation}
The metric evolution equation \eqref{eq:metric-evolution} therefore
gives
\begin{equation}\label{eq:prop39m-metric-comparison}
|\partial_sg_m|_{g_m}\leq C\sigma_m^2,
\qquad
e^{-C\sigma_m^2}g_m(1)
\leq g_m(s)\leq e^{C\sigma_m^2}g_m(1).
\end{equation}

\smallskip
\noindent\emph{Step 2. A common graphical representation.}
Let $y_m:=F_m(p_m,1)$ and choose a $J$-adapted
$\bar g$-orthonormal frame at $y_m$ whose first two vectors span
$d(F_m(\cdot,1))_{p_m}(T_{p_m}\Sigma)$. Define
\[
E_m(y):=\exp_{y_m}^{\bar g}(\sqrt{d_m}\,y),
\qquad
\widetilde g_m:=E_m^*\bar g_m,
\qquad
J_m:=E_m^*J,
\qquad
\widetilde F_m:=E_m^{-1}\circ F_m.
\]
Here $B_R^{\mathbb R^k}(0)$ denotes the Euclidean ball of radius
$R$ in $\mathbb R^k$. Since $d_m\leq1/m$, we have
\[
d_m^{-1/2}\operatorname{inj}_{\bar g}(M)\longrightarrow\infty
\qquad\text{as }m\to\infty.
\]
Hence these rescaled normal-coordinate neighborhoods exhaust
$\mathbb R^4$. On every fixed ball and
for every integer $q\geq0$,
\begin{equation}\label{eq:prop39m-normal-coordinates}
\|\widetilde g_m-g_{\mathrm{Euc}}\|_{C^q(B_R^{\mathbb R^4}(0))}
\leq C_{q,R}d_m.
\end{equation}
At $s=1$, we have $\widetilde F_m(p_m,1)=0$, and the tangent
plane of $\widetilde F_m(\cdot,1)$ at $p_m$ is
$\mathbb R^2\times\{0\}$. Moreover,
$\widetilde g_m(0)=g_{\mathrm{Euc}}$ and
$\Gamma(\widetilde g_m)(0)=0$.

Let $\operatorname{pr}_{\mathrm h}:\mathbb R^2\times\mathbb R^2
\to\mathbb R^2$ be the horizontal projection. The velocity bound in
\eqref{eq:prop39m-rescaled-bounds} gives
\begin{equation}\label{eq:prop39m-displacement}
|\widetilde F_m(q,s)-\widetilde F_m(q,1)|\leq C\sigma_m
\end{equation}
whenever the two points lie in a fixed coordinate ball. Set
$r_m:=\sqrt{\sigma_m}$. Applying
Lemma~\ref{lem:chord-tangent} to a $g_m(s)$-geodesic of length
$r_m$ issuing from $p_m$ and to a minimizing $g_m(1)$-geodesic joining the same
two points of $\Sigma$, and using
\eqref{eq:prop39m-metric-comparison} and
\eqref{eq:prop39m-displacement} at both endpoints, yields
\begin{equation}\label{eq:prop39m-base-angle}
\sup_{0\leq s\leq1}
\angle\!\left(
d(\widetilde F_m(\cdot,s))_{p_m}(T_{p_m}\Sigma),
\mathbb R^2\times\{0\}
\right)
\leq
C\left[
r_m(\sigma_m+d_m)+\frac{\sigma_m}{r_m}+\sigma_m^2
\right]
=o(1).
\end{equation}

Fix $R>0$ and set
$z_m(s):=\operatorname{pr}_{\mathrm h}\widetilde F_m(p_m,s)$.
For $x\in B_R^{\mathbb R^2}(0)$, lift the segment from $z_m(s)$ to
$x$ by solving
\begin{equation}\label{eq:prop39m-radial-lift}
q'(t)=
\left[
d\bigl(\operatorname{pr}_{\mathrm h}
       \widetilde F_m(\cdot,s)\bigr)_{q(t)}
\right]^{-1}
(x-z_m(s)),
\qquad q(0)=p_m.
\end{equation}
The inverse exists initially by \eqref{eq:prop39m-base-angle}. As
long as the tangent-plane angle is below $1/4$, the lifted curve has
length at most $C_R$, and the Gauss formula changes its tangent plane
by at most $C_R(\sigma_m+d_m)$. Hence
\eqref{eq:prop39m-base-angle} and a first-exit argument prevent both
loss of invertibility and exit from a fixed coordinate ball for all
large $m$. The lift therefore reaches its endpoint. Uniqueness and smooth dependence select the local branch of the
immersion determined by $p_m$, even at a
self-intersection. Moreover, if $R_1<R_2$, then the construction
on $B_{R_2}(0)$ restricts to that on $B_{R_1}(0)$.

Writing the selected sheet as
\[
\widetilde F_m(q_{m,R}(x,s),s)=(x,u_{m,R}(x,s)),
\]
we have $u_{m,R},Du_{m,R}\to0$ uniformly on
$B_R^{\mathbb R^2}(0)\times[0,1]$. We regard the Euclidean spatial
Hessian $D^2u_{m,R}$ as taking values in
\[
\mathscr E
:=\mathbb R^2\otimes
\operatorname{Sym}^2\bigl((\mathbb R^2)^*\bigr),
\]
equipped with its Frobenius norm. Suppressing the subscript $R$
in the following local computation, the graph Gauss formula gives
\begin{equation}\label{eq:prop39m-graph-gauss}
(h_m)_{ij}
=\pi_{\widetilde g_m}^{\perp}
\left[
(0,\partial_i\partial_j u_m)
+\Gamma(\widetilde g_m)
 \bigl(\partial_i(x,u_m),\partial_j(x,u_m)\bigr)
\right].
\end{equation}
Since $Du_m=o(1)$, this implies
$|D^2u_m|\leq(1+o(1))|h_m|+C_Rd_m$. At $(0,1)$ the
coordinates are normal and $Du_m(0,1)=0$, so the Hessian array is
exactly the second fundamental form array. A diagonal choice now
gives $R_m\to\infty$ and smooth graphs
\[
u_m:B_{R_m}^{\mathbb R^2}(0)\times[0,1]\to\mathbb R^2
\]
such that
\begin{equation}\label{eq:prop39m-graph-bounds}
\begin{gathered}
|D^2u_m|\leq2\sqrt L\,\sigma_m,
\qquad
|D^2u_m(\cdot,0)|\leq2\sigma_m,\\
u_m(0,1)=0,
\qquad
Du_m(0,1)=0,
\qquad
|D^2u_m(0,1)|=\sqrt L\,\sigma_m.
\end{gathered}
\end{equation}

\smallskip
\noindent\emph{Step 3. The normalized graph system.}
Let $\widetilde\ff_m$ be the coordinate representation of $\ff_m$ and,
with respect to
$\mathbb R^4=\mathbb R^2_{\mathrm h}\oplus
\mathbb R^2_{\mathrm v}$, write
\[
\widetilde\ff_m
=\bigl((\widetilde\ff_m)^{\mathrm h},
       (\widetilde\ff_m)^{\mathrm v}\bigr),
\qquad
(\widetilde\ff_m)^{\mathrm h}
=\bigl((\widetilde\ff_m)^{\mathrm h}_1,
       (\widetilde\ff_m)^{\mathrm h}_2\bigr).
\]
All components are evaluated along $(x,u_m(x,s))$. Keeping the
horizontal graph coordinate fixed adds a tangential
reparametrization, and hence
\begin{equation}\label{eq:prop39m-horizontal-gauge}
\partial_su_m
=(\widetilde\ff_m)^{\mathrm v}
-(\widetilde\ff_m)^{\mathrm h}_i\partial_i u_m.
\end{equation}
Substituting \eqref{eq:prop39m-graph-gauss} into this identity and
collecting all $D^2u_m$ terms into the principal part gives
\begin{equation}\label{eq:prop39m-graph-system}
\partial_su_m
=\widetilde{\mathcal A}_m^{ij}(x,s)\partial_i\partial_j u_m
+d_m\widetilde{\mathcal B}_m(x,s).
\end{equation}
The normal-coordinate estimate
\eqref{eq:prop39m-normal-coordinates} and $\bar\nabla J=0$ imply that,
on every fixed ball,
\[
\|\Gamma(\widetilde g_m)\|_{C^1(B_R^{\mathbb R^4}(0))}
+\|\partial J_m\|_{C^1(B_R^{\mathbb R^4}(0))}
\leq C_Rd_m.
\]
Together with \eqref{eq:prop39m-graph-bounds}, these estimates
imply that, for every $\gamma\in(0,1)$,
$\widetilde{\mathcal A}_m^{ij}$ and
$\widetilde{\mathcal B}_m$ have uniform spatial $C^\gamma$ bounds
on smaller fixed cylinders, as required for the interior Schauder
estimate used below.

Let $P_m^\perp$ denote the principal symbol of the normal
linearization at $(x,s)$. The calculation preceding
Definition~\ref{def:heat_kernel} gives, in adapted orthonormal
frames,
\[
\left\langle
P_m^\perp(x,s,\xi)\eta,\eta
\right\rangle
=
a_m|\xi|^2|\eta|^2
+
(1-a_m)(\xi_2\eta_3+\xi_1\eta_4)^2,
\qquad
a_m:=
\frac{\cos^2\alpha_m}
{\cos^2\alpha_m+\beta\sin^2\alpha_m}.
\]
The principal symbol
\[
\widetilde P_m(x,s,\xi)
:=
\widetilde{\mathcal A}_m^{ij}(x,s)\xi_i\xi_j
\]
of the graph system is obtained from $P_m^\perp$ by the
zeroth-order changes from the graph coordinate frames to adapted
orthonormal tangent and normal frames.
Since $Du_m\to0$ and $\widetilde g_m\to g_{\mathrm{Euc}}$
locally uniformly, orthonormalizing the graph tangent frame
and the normal projections of the vertical frame gives changes
converging to the identity. The subsequent changes to adapted
orthonormal frames are orthogonal and preserve the symbol bounds.
Hence the graph systems are uniformly parabolic on each fixed
compact subcylinder for all sufficiently large $m$.

After passing to a
subsequence, $J_m(0)$ converges, and
$\widetilde g_m\to g_{\mathrm{Euc}}$,
$J_m-J_m(0)\to0$, and $u_m,Du_m\to0$ show that
$\widetilde{\mathcal A}_m^{ij}$ converges locally uniformly to constant matrices
$\widetilde{\mathcal A}_\infty^{ij}$. Their limiting symbol is symmetric and
satisfies
\begin{equation}\label{eq:prop39m-limit-symbol}
\kappa|\xi|^2I
\leq\widetilde{\mathcal A}_\infty^{ij}\xi_i\xi_j
\leq|\xi|^2I.
\end{equation}

Set
\[
v_m:=\frac{u_m}{\sigma_m},
\qquad
\mathcal R_m:=\frac{d_m}{\sigma_m}\widetilde{\mathcal B}_m.
\]
Then
\begin{equation}\label{eq:prop39m-normalized-system}
\partial_sv_m
=\widetilde{\mathcal A}_m^{ij}\partial_i\partial_jv_m+\mathcal R_m,
\qquad
\|\mathcal R_m\|_{C_x^\gamma(K\times[0,1])}\to0
\end{equation}
for every compact $K\subset\mathbb R^2$. The velocity bound and
\eqref{eq:prop39m-horizontal-gauge} give a uniform local bound for
$\partial_sv_m$, while \eqref{eq:prop39m-graph-bounds} becomes
\begin{equation}\label{eq:prop39m-normalized-bounds}
\begin{gathered}
|D^2v_m|\leq2\sqrt L,
\qquad
|D^2v_m(\cdot,0)|\leq2,\\
v_m(0,1)=0,
\qquad
Dv_m(0,1)=0,
\qquad
|D^2v_m(0,1)|=\sqrt L=4C_0.
\end{gathered}
\end{equation}

\smallskip
\noindent\emph{Step 4. Parabolic compactness up to $s=1$.}
Fix $R>0$ and $\varepsilon>0$. Taylor's formula at $s=1$ and the
bound for $\partial_sv_m$ give a uniform $L^2$ bound for $v_m$ on
$B_{R+2}^{\mathbb R^2}(0)\times[0,1]$. The interior Schauder estimate for parabolic systems
\cite[Theorem~2.2(I)]{DongZhang}, applied locally and
followed by a finite covering argument, gives a uniform
$C^{\gamma/2,\gamma}$ bound for $D^2v_m$ on
$B_R^{\mathbb R^2}(0)\times[\varepsilon,1)$.
By continuity, the same bound extends to $s=1$. Together with the standard finite-difference estimate for $Dv_m$,
this yields, after a diagonal extraction,
\[
v_m\longrightarrow v_\infty
\quad\text{locally on }\mathbb R^2\times(0,1],
\]
with locally uniform convergence of the spatial derivatives through
order two. Since
\[
\widetilde{\mathcal A}_m^{ij}
\longrightarrow
\widetilde{\mathcal A}_\infty^{ij},
\qquad
\mathcal R_m\longrightarrow0
\]
locally uniformly, passing to the limit in
\eqref{eq:prop39m-normalized-system} gives
\begin{equation}\label{eq:prop39m-limit-system}
\partial_sv_\infty=L_\infty v_\infty,
\qquad
L_\infty
:=
\widetilde{\mathcal A}_\infty^{ij}
\partial_i\partial_j,
\qquad
|D^2v_\infty(0,1)|=4C_0.
\end{equation}

\smallskip
\noindent\emph{Step 5. The initial trace and the contradiction.}
By the Banach--Alaoglu theorem and a diagonal argument on the
expanding balls, there exist a subsequence, still denoted by
$v_m$, and
\[
W_0\in L^\infty(\mathbb R^2;\mathscr E)
\]
such that
\begin{equation}\label{eq:prop39m-initial-trace}
D^2v_m(\cdot,0)\stackrel{*}{\rightharpoonup}W_0
\quad\text{in }L^\infty_{\mathrm{loc}}
(\mathbb R^2;\mathscr E),
\qquad
\|W_0\|_{L^\infty}\leq2.
\end{equation}
For every compactly supported smooth $\mathscr E$-valued test field
$\varphi$, two integrations by parts and the local bound for
$\partial_sv_m$ give
\[
\left|
\left\langle
D^2v_m(\cdot,s)-D^2v_m(\cdot,0),\varphi
\right\rangle
\right|
\leq C_\varphi s.
\]
Letting first $m\to\infty$ and then $s\downarrow0$, we obtain
\[
D^2v_\infty(\cdot,s)\longrightarrow W_0
\quad\text{in }\mathcal D'(\mathbb R^2;\mathscr E).
\]

Set $W:=D^2v_\infty$. The bound in
\eqref{eq:prop39m-normalized-bounds} passes to $W$, and positive-time
interior regularity makes $W$ a bounded classical solution. Since
$L_\infty$ has constant coefficients,
\[
\partial_sW
=\bigl(L_\infty\otimes
I_{\operatorname{Sym}^2((\mathbb R^2)^*)}\bigr)W.
\]
Let $G_s^\infty$ be the matrix heat kernel of $L_\infty$.
Corollary~\ref{cor:frozen-semigroup-bound}, applied to the first
factor, yields
\[
W(\cdot,s)
=\left(
G_s^\infty\otimes
I_{\operatorname{Sym}^2((\mathbb R^2)^*)}
\right)*W_0.
\]
Tensoring with the identity preserves the operator norm, so
Lemma~\ref{lem:frozen-kernel-bound} and
\eqref{eq:prop39m-initial-trace} imply
\[
\|D^2v_\infty(\cdot,1)\|_{L^\infty}
\leq C_0\|W_0\|_{L^\infty}
\leq2C_0.
\]
This contradicts \eqref{eq:prop39m-limit-system}. Hence the required
$\theta>0$ exists. 
\end{proof}

\begin{proof}[Proof of Theorem~\ref{theorem:blowup_lower_rate}]
Let $L$ and $\theta$ be the constants given by
Proposition~\ref{prop:curvature-persistence}.
Fix $t_0\in[0,T)$ and set $K:=1+Q(t_0)$.
If $T-t_0<\theta/K$, then
Proposition~\ref{prop:curvature-persistence} gives
\[
Q(t)<LK
\qquad (t_0\leq t<T).
\]
Since the solution is smooth on $[0,t_0]$, it follows that
$\sup_{0\leq t<T}Q(t)<\infty$.
By Theorem~\ref{thm:extension}, the solution then extends smoothly
beyond $T$, contradicting the maximality of $T$.

Thus $T-t_0\geq\theta/K$. If $0<T-t_0\leq\theta/2$, then
\[
Q(t_0)
\geq\frac{\theta}{T-t_0}-1
\geq\frac{\theta}{2(T-t_0)}.
\]
Since $t_0$ was arbitrary subject to
$0<T-t_0\leq\theta/2$, this proves
\eqref{eq:blowup_lower_rate}.
\end{proof}

We also record a Euclidean version of
Proposition~\ref{prop:curvature-persistence} for complete
immersions, which will be used in Section~\ref{sec:type-I}.

\begin{lemma}[Euclidean curvature persistence]
\label{lem:euclidean-persistence}
Fix $\beta>0$ and $0<v_*\leq1$. There exist constants
$L_E>1$ and $\theta_E>0$, depending only on $\beta$ and $v_*$,
with the following property. Let $S$ be a possibly noncompact
smooth oriented surface without boundary, and let
$F:S\times I\to\R^4$ be a smooth solution of
\eqref{eq:gradient-flow} for a constant orthogonal complex
structure, where $I$ is an open interval containing $[0,S_0]$
and $S_0>0$. Assume that every induced metric is complete,
$\cos\alpha\geq v_*$, and
\[
\sup_{S\times I'}|\nabla^m h|<\infty
\qquad
(m=0,1,2,\ I'\Subset I).
\]
For $q>0$, set $I_q:=[0,S_0]\cap[0,\theta_E/q)$.
If $\sup_S|h(\cdot,0)|^2\leq q$, then
\begin{equation}\label{eq:euclidean-persistence}
\sup_{S\times I_q}|h|^2\leq L_Eq.
\end{equation}
\end{lemma}

\begin{proof}
The proof of Proposition~\ref{prop:curvature-persistence}
applies with the following modifications.
The evolution equation for $h$ and the velocity formula show
that the assumed derivative bounds give a uniform bound for
$\partial_t|h|^2$ on $S\times I'$ for every $I'\Subset I$.
Thus $t\mapsto\sup_S|h(\cdot,t)|^2$ is locally Lipschitz,
so the first hitting time construction remains valid. At the first hitting time $b_m$, choose $p_m\in S$ such that
\[
|h(p_m,b_m)|^2
\geq (1-1/m)\sup_S |h(\cdot,b_m)|^2.
\]
Completeness suffices to continue the bounded-length
radial lifts.

By the Euclidean parabolic rescaling
$\widehat F(x,s):=\sqrt q\,F(x,s/q)$, it suffices to treat $q=1$.
The ambient lower-order term in \eqref{eq:prop39m-graph-system}
vanishes, and the rescaling, compactness, and heat-kernel arguments
of Proposition~\ref{prop:curvature-persistence} apply as above.
\end{proof}

\section{Exclusion of Type I singularities}
\label{sec:type-I}
\label{sec:special-type-I}

To prove Theorem~\ref{thm:main}, we combine the curvature estimates
of Section~\ref{sec:proof} with the weighted monotonicity formula
of Han, Li, and Sun~\cite{HLSflow}.

\begin{definition}\label{def:type-I}
A finite singular time $T$ of the flow \eqref{eq:gradient-flow}
is called \emph{Type~I} if there is a constant $C_I>0$ such that
\begin{equation}\label{eq:typeI-bound}
\max_{x\in\Sigma}|h(x,t)|^2\leq\frac{C_I}{T-t}
\qquad (0\leq t<T).
\end{equation}
\end{definition}

\begin{theorem}\label{thm:no-type-I}
Let $(M,\omega,J)$ be a compact K\"ahler surface,
let $\Sigma$ be closed and oriented, and fix $\beta>0$.
Suppose that $F:\Sigma\times[0,T)\to M$ is a smooth solution
of \eqref{eq:gradient-flow}, with $0<T<\infty$ and
$\min_\Sigma\cos\alpha(\cdot,0)>0$. If \eqref{eq:typeI-bound}
holds, then
\begin{equation}\label{eq:typeI-bounded}
 \sup_{0\leq t<T}\max_\Sigma|h(\cdot,t)|^2<\infty,
\end{equation}
and the solution extends smoothly across $T$. Consequently,
the flow develops no Type~I singularities.
\end{theorem}

We work under the assumptions of Theorem~\ref{thm:no-type-I}.
Since $T<\infty$, estimate \eqref{eq:2.3} gives a constant
$v_*>0$ such that
\begin{equation}\label{eq:typeI-angle}
v:=\cos\alpha\geq v_*
\qquad\text{on }\Sigma\times[0,T).
\end{equation}
By \eqref{eq:gradient-flow-f}, the normal velocity $\ff$ is a
convex combination of $\HH$ and $-\VV$. Using
$|\HH|,|\VV|\leq\sqrt{2}\,|h|$, we obtain
\begin{equation}\label{eq:typeI-speed}
|\ff|^2\leq 2|h|^2,
\qquad
|\ff|^2\leq|\HH|^2+|\VV|^2.
\end{equation}

Combining \eqref{eq:typeI-bound} with \eqref{eq:typeI-speed} gives
$|\partial_t F(x,t)|\leq\sqrt{2C_I}(T-t)^{-1/2}$.
Hence, for $x\in\Sigma$ and $0\leq t<t'<T$,
\[
d_{\bar g}\bigl(F(x,t),F(x,t')\bigr)
\leq \sqrt{2C_I}\int_t^{t'}(T-u)^{-1/2}\,du
=2\sqrt{2C_I}\bigl(\sqrt{T-t}-\sqrt{T-t'}\bigr).
\]
Since this estimate is uniform in $x$ and $M$ is compact,
$F(\cdot,t)$ converges uniformly as $t\uparrow T$ to a
continuous map $F_T:\Sigma\to M$. Letting $t'\uparrow T$ gives
\begin{equation}\label{eq:typeI-terminal}
d_{\bar g}\bigl(F(x,t),F_T(x)\bigr)
\leq 2\sqrt{2C_I}\sqrt{T-t}.
\end{equation}

By the calculation in Section~\ref{sec:proof} and the angle bound
\eqref{eq:typeI-angle}, the positive principal symbol $P$ of the
normal linearization of \eqref{eq:gradient-flow} satisfies
\begin{equation}\label{eq:typeI-symbol}
\kappa|\xi|^2 I\leq P(\xi)\leq |\xi|^2 I,
\qquad
\kappa:=\frac{v_*^2}{v_*^2+\beta(1-v_*^2)}>0.
\end{equation}

To study a possible singularity at $T$, let $\tau_j\downarrow0$
and define the parabolic rescalings
\[
\bar g_j:=\tau_j^{-1}\bar g,
\qquad
F_j(x,s):=F\bigl(x,T-\tau_j+\tau_j s\bigr).
\]
Let $h_j$ denote the second fundamental form of $F_j$ with respect
to $\bar g_j$. Since parabolic rescaling preserves the K\"ahler
angle, \eqref{eq:typeI-symbol} holds for each $F_j$ with the
same constant $\kappa$.

\begin{lemma}\label{lem:typeI-compactness}
Fix $0<b<1$, and let $x_j\in\Sigma$ and $X_j\in M$ satisfy
\[
\sup_j d_{\bar g_j}\bigl(F_j(x_j,b),X_j\bigr)<\infty.
\]
In rescaled normal coordinates centered at $X_j$, a subsequence of
$(F_j,x_j)$ converges smoothly in the pointed sense to a complete
$L_\beta$-gradient flow
$F_\infty:S_\infty\times I\to(\R^4,J_\infty)$, where $I$ is an
open interval containing $[0,b]$ and $J_\infty$ is a constant
orthogonal complex structure. The limit satisfies
$\cos\alpha_\infty\geq v_*$ and
\begin{equation}\label{eq:typeI-derivative-bounds}
\sup_{S_\infty\times I'}|\nabla^m h_\infty|<\infty
\qquad (m\geq0,\ I'\Subset I).
\end{equation}
\end{lemma}

\begin{proof}
Choose an open interval $I$ with $[0,b]\subset I\Subset(-\infty,1)$,
and write $g_j(s)=F_j(\cdot,s)^*\bar g_j$. For all large $j$,
the Type~I bound gives $|h_j|^2\leq C_I/(1-s)$ on $\Sigma\times I$.
By \eqref{eq:typeI-speed} and \eqref{eq:metric-evolution},
$g_j(s)$ is uniformly equivalent to $g_j(b)$ on $I$.

The curvature and speed bounds give local graphs of uniformly small
slope on cylinders of uniform size, as in Step~2 of the proof of
Proposition~\ref{prop:curvature-persistence}. The graph Gauss
formula bounds $D^2u_j$, so the coefficients and the lower-order
term are uniformly Lipschitz in space. The ellipticity
bound \eqref{eq:typeI-symbol} implies uniform parabolicity of the
graph systems. Applying the interior estimate of Step~4 and
iterating after spatial differentiation gives
\[
\sup_{\Sigma\times I'}|\nabla^m h_j|\leq C_{m,I'}
\qquad (m\geq0,\ I'\Subset I).
\]
The velocity formula then bounds all spatial covariant derivatives
of $\partial_s F_j$ on the same intervals.

In rescaled normal coordinates centered at $X_j$, the ambient
metrics and complex structures subconverge smoothly to the Euclidean
metric and a constant orthogonal complex structure. By hypothesis,
$p_j:=F_j(x_j,b)$ stays bounded in these coordinates. After passing
to a subsequence, $p_j$ converges. With ambient basepoints $p_j$,
Smith's compactness argument \cite[Theorem~1.2]{Smith}, applied at
all orders by diagonal extraction in the graph charts, gives smooth
pointed convergence at time $b$ to a complete immersion.

Let $\Phi_j:U_j\to\Sigma$ denote the convergence maps at time $b$,
where the domains $U_j$ exhaust $S_\infty$ and
$\Phi_j(x_\infty)=x_j$ for the limit basepoint $x_\infty$.
Following the argument in the proof of \cite[Theorem~2.3]{ChenHe},
we keep these maps fixed in time. On each compact time interval
containing $b$, the speed bound keeps the images of compact subsets
of $S_\infty$ in bounded coordinate regions. Starting at $s=b$,
we use the metric and connection variation formulas together with
the curvature and velocity estimates to bound all spatial derivatives
of $F_j(\Phi_j(\cdot),s)$ in the fixed charts. The flow equation
then gives the time and mixed derivative bounds. By the
Arzel\`a--Ascoli theorem and a diagonal argument on a compact
exhaustion, a subsequence converges smoothly on compact subsets
of $S_\infty\times I$. The flow equation and the angle bound pass to the limit.
By smooth local convergence, the same bound $C_{m,I'}$ holds
at every point of $S_\infty\times I'$, which gives
\eqref{eq:typeI-derivative-bounds}. Since $g_\infty(s)$ is
uniformly equivalent to the complete metric $g_\infty(b)$,
it is complete for every $s\in I$.
\end{proof}

An initial curvature bound on expanding extrinsic balls gives
a global bound on the limit surface. Suppose that, for some
$\varepsilon\geq0$ and some sequence $R_j\to\infty$,
\begin{equation}\label{eq:typeI-expanding-initial}
|h_j(x,0)|^2\leq\varepsilon
\qquad\text{whenever }d_{\bar g_j}(F_j(x,0),X_j)\leq R_j.
\end{equation}
Fix $z\in S_\infty$. By convergence at $s=0$ in the rescaled
normal coordinates centered at $X_j$, the distances
$d_{\bar g_j}(F_j(\Phi_j(z),0),X_j)$ remain bounded.
Since $R_j\to\infty$, \eqref{eq:typeI-expanding-initial} applies
at $\Phi_j(z)$ for all sufficiently large $j$. Passing to the
limit gives $|h_\infty(z,0)|^2\leq\varepsilon$, and therefore
\begin{equation}\label{eq:typeI-limit-initial}
\sup_{S_\infty}|h_\infty(\cdot,0)|^2\leq\varepsilon.
\end{equation}

We combine the initial bound \eqref{eq:typeI-limit-initial} with
Lemma~\ref{lem:euclidean-persistence} to propagate a curvature
bound to a smaller extrinsic ball at a later time.

\begin{lemma}\label{lem:typeI-finite-step}
Let $L_E$ and $\theta_E$ be the constants in
Lemma~\ref{lem:euclidean-persistence}, and fix
$0<\varepsilon<\theta_E/2$ and $K>4L_E$.
There exist $R>1$ and $\tau_*\in(0,T)$ with the following property.
For $X\in M$ and $0<\tau<\tau_*$, set
$t_0=T-\tau$ and $t_1=T-\tau/K$. If
\begin{equation}\label{eq:typeI-initial-small}
|h(x,t_0)|^2\leq\frac{\varepsilon}{\tau}
\qquad\text{whenever }d_{\bar g}(F(x,t_0),X)\leq R\sqrt{\tau},
\end{equation}
then
\begin{equation}\label{eq:typeI-step-conclusion}
|h(x,t_1)|^2\leq\frac{2L_E\varepsilon}{\tau}
\qquad\text{whenever }d_{\bar g}(F(x,t_1),X)\leq\sqrt{\tau}.
\end{equation}
In particular, at each such point,
\begin{equation}\label{eq:typeI-normalized-gain}
(T-t_1)|h(x,t_1)|^2
\leq\frac{2L_E}{K}\varepsilon<\frac{\varepsilon}{2}.
\end{equation}
\end{lemma}

\begin{proof}
Suppose that no such $R$ and $\tau_*$ give
\eqref{eq:typeI-step-conclusion}. Then there are counterexamples
with $R_j\to\infty$, $\tau_j\downarrow0$, and centers $X_j\in M$.
Rescale by $\tau_j$ as above and set $b=1-K^{-1}$, so that
$s=0$ and $s=b$ correspond to $T-\tau_j$ and $T-\tau_j/K$,
respectively. The initial estimate \eqref{eq:typeI-initial-small}
becomes \eqref{eq:typeI-expanding-initial}, while the failure of
\eqref{eq:typeI-step-conclusion} gives points $x_j\in\Sigma$ with
\[
d_{\bar g_j}(F_j(x_j,b),X_j)\leq1,
\qquad
|h_j(x_j,b)|^2>2L_E\varepsilon.
\]

Lemma~\ref{lem:typeI-compactness} gives a complete pointed limit.
By \eqref{eq:typeI-limit-initial} and smooth convergence at the
marked point,
\[
\sup_{S_\infty}|h_\infty(\cdot,0)|^2\leq\varepsilon,
\qquad
|h_\infty(x_\infty,b)|^2\geq2L_E\varepsilon.
\]
The angle and derivative bounds on the limit allow us to apply
Lemma~\ref{lem:euclidean-persistence} with $q=\varepsilon$ and
$S_0=b$. Since $b<1<\theta_E/\varepsilon$, it follows that
\[
\sup_{S_\infty}|h_\infty(\cdot,b)|^2\leq L_E\varepsilon,
\]
contradicting the lower bound at $x_\infty$. This proves
\eqref{eq:typeI-step-conclusion}. Multiplying by $T-t_1=\tau/K$
and using $K>4L_E$ gives \eqref{eq:typeI-normalized-gain}.
\end{proof}

Lemma~\ref{lem:typeI-finite-step} propagates a local curvature bound
between two time slices. We next use the weighted monotonicity formula
of Han, Li, and Sun to prove uniform decay of $(T-t)|h|^2$ on balls
of radius $B\sqrt{T-t}$ about a fixed point, for each fixed $B>0$.

Fix $X\in M$ and $0<2r<\operatorname{inj}_{\bar g}(M)$. Choose
$\phi\in C_c^\infty(B_{2r}(X))$ with $0\leq\phi\leq1$ and
$\phi\equiv1$ on $B_r(X)$. In normal coordinates $y$ centered at
$X$, define the backward Gaussian weight and the weighted integral by
\[
\rho_X(y,t)
:=\frac{1}{4\pi(T-t)}
  \exp\!\left(-\frac{|y|^2}{4(T-t)}\right),
\qquad
\Psi_X(t):=\int_{\Sigma_t}v^{-q}\phi\rho_X\dm_t,
\]
where $0\leq t<T$ and $q>0$ is fixed sufficiently large.
Integrals over $\Sigma_t$ are taken on the parameter surface
$\Sigma$ with the area form induced by $F_t$.

Dropping one nonnegative square from the localized estimate in
the proof of \cite[Theorem~3.3]{HLSflow}, we obtain
\begin{equation}\label{eq:typeI-monotonicity}
\begin{aligned}
\frac{d}{dt}\bigl(e^{c_1\sqrt{T-t}}\Psi_X(t)\bigr)
\leq{}& -e^{c_1\sqrt{T-t}}
\int_{\Sigma_t} v^{-q}\phi\rho_X
\left[
  \frac14\left|\frac{y^\perp}{T-t}+\ff+\HH\right|^2
  +\frac78|\HH|^2+\frac78|\VV|^2
\right]\dm_t\\
& +c_2e^{c_1\sqrt{T-t}}.
\end{aligned}
\end{equation}
Here $c_1,c_2>0$ are independent of $t$, and $y^\perp$ denotes
the normal component of the radial vector field along $F_t$.

\begin{lemma}\label{lem:typeI-cone-decay}
Fix $X\in M$ and $B>0$, and set
\[
N_X(t):=\{x\in\Sigma:
d_{\bar g}(F(x,t),X)\leq B\sqrt{T-t}\}.
\]
Then
\begin{equation}\label{eq:typeI-cone-decay}
\lim_{t\uparrow T}(T-t)\sup_{x\in N_X(t)}|h(x,t)|^2=0,
\end{equation}
where the supremum is taken to be zero if $N_X(t)$ is empty.
\end{lemma}

\begin{proof}
If \eqref{eq:typeI-cone-decay} fails, there exist $\delta>0$,
times $t_j\uparrow T$, and points $x_j\in\Sigma$ such that
\begin{equation}\label{eq:typeI-cone-failure}
d_{\bar g}(F(x_j,t_j),X)\leq B\sqrt{T-t_j},
\qquad
(T-t_j)|h(x_j,t_j)|^2\geq\delta.
\end{equation}
Set $\tau_j=2(T-t_j)$ and $b=1/2$. In the parabolic rescalings
defined above, time $s=b$ corresponds to $t_j$. Hence
\[
d_{\bar g_j}(F_j(x_j,b),X)\leq B/\sqrt{2},
\qquad
|h_j(x_j,b)|^2\geq2\delta.
\]
Lemma~\ref{lem:typeI-compactness}, applied with $X_j=X$, gives
a smooth pointed limit in rescaled normal coordinates centered
at $X$, with $|h_\infty(x_\infty,b)|^2\geq2\delta$.

Integrating \eqref{eq:typeI-monotonicity} in time and using
$\Psi_X\geq0$, we see that the nonnegative dissipation term is
integrable up to $T$, since the error term is integrable there.
Its integral over $[T-\tau_j,t_j]$ therefore tends to zero.
Under the rescaling $t=T-\tau_j+\tau_j s$, this interval
corresponds to $[0,b]$. Writing $F_j$ in the rescaled normal
coordinates, the Gaussian weight along $F_j$ is
\[
\rho_j(x,s)
:=\tau_j\rho_X\bigl(\sqrt{\tau_j}\,F_j(x,s),t\bigr)
=\frac{1}{4\pi(1-s)}
  \exp\!\left(-\frac{|F_j(x,s)|^2}{4(1-s)}\right).
\]
On each compact parameter patch $K\Subset S_\infty$ and for
$0\leq s\leq b$, smooth convergence gives a uniform positive
lower bound for $\rho_j$ on $\Phi_j(K)$, and the rescaled cutoff
equals one there for all large $j$. 
Since
\[
|F_j^\perp|
\leq (1-s)\left(
\left|\frac{F_j^\perp}{1-s}+\ff_j+\HH_j\right|
+|\ff_j|+|\HH_j|
\right),
\]
the speed estimate \eqref{eq:typeI-speed} shows that the rescaled
square bracket in \eqref{eq:typeI-monotonicity} controls
$|\HH_j|^2+|F_j^\perp|^2$ on $[0,b]$. Restricting the
rescaled integral to $\Phi_j(K)$ gives
\[
\int_0^b\int_{\Phi_j(K)}
\bigl(|\HH_j|^2+|F_j^\perp|^2\bigr)\,d\mu_s^j\,ds
\longrightarrow0.
\]
Smooth convergence and continuity up to $s=b$ now yield
\begin{equation}\label{eq:typeI-cone-identities}
\HH_\infty=0,\qquad F_\infty^\perp=0
\quad\text{on }S_\infty\times[0,b].
\end{equation}

At time $b$, wherever $F_\infty\ne0$, write
$F_\infty=dF_\infty(W)$ with $W\ne0$. Differentiating in a
tangent direction $Z$ and taking normal components gives
$h_\infty(Z,W)=0$. Since the surface is two-dimensional and
$\HH_\infty=0$, it follows that $h_\infty=0$ there. The zero
set of the position vector has empty interior because
$F_\infty(\cdot,b)$ is an immersion. By continuity, $h_\infty$
vanishes on the entire slice, contradicting
$|h_\infty(x_\infty,b)|^2\geq2\delta$.
\end{proof}

To apply Lemma~\ref{lem:typeI-cone-decay}, we use local curvature
propagation to select high-curvature points whose images approach
a single center at the rate $O(\sqrt{T-t})$.

\begin{lemma}\label{lem:typeI-chain}
Suppose that $T<\infty$ is the maximal time of smooth existence
and that \eqref{eq:typeI-bound} holds. Then there exist constants
$\varepsilon,B>0$, a point $X_\infty\in F_T(\Sigma)$, and
sequences $t_n\uparrow T$ and $x_n\in\Sigma$ such that
\begin{equation}\label{eq:typeI-essential-chain}
(T-t_n)|h(x_n,t_n)|^2\geq\varepsilon,
\qquad
d_{\bar g}(F(x_n,t_n),X_\infty)\leq B\sqrt{T-t_n}.
\end{equation}
\end{lemma}

\begin{proof}
Let $\theta$ be the constant in
Theorem~\ref{theorem:blowup_lower_rate}. Choose
$0<\varepsilon<\frac12\min\{\theta,\theta_E\}$ and $K>4L_E$.
Let $R,\tau_*$ be given by Lemma~\ref{lem:typeI-finite-step},
and choose $0<\tau_0<\min\{\tau_*,\theta/2\}$. For $n\geq0$, set
\begin{equation}\label{eq:typeI-concentration-sets}
\tau_n=\tau_0K^{-n},\qquad t_n=T-\tau_n,\qquad
E_n=\{x\in\Sigma:\tau_n|h(x,t_n)|^2\geq\varepsilon\}.
\end{equation}
Each $E_n$ is compact and is nonempty by
\eqref{eq:blowup_lower_rate}, since
$\tau_n Q(t_n)\geq\theta/2>\varepsilon$.

For each $x\in E_{n+1}$, the contrapositive of
Lemma~\ref{lem:typeI-finite-step}, with $\tau=\tau_n$ and
$X=F(x,t_{n+1})$, gives a predecessor $z\in E_n$ satisfying
\begin{equation}\label{eq:typeI-chain-jump}
d_{\bar g}(F(z,t_n),F(x,t_{n+1}))\leq R\sqrt{\tau_n}.
\end{equation}
For each $N$, start with a point of $E_N$ and choose successive
predecessors down to $E_0$. By compactness of each $E_n$, a
diagonal subsequence of these finite chains converges at every
fixed level $n$ to a point $x_n\in E_n$. Passing to the limit
in \eqref{eq:typeI-chain-jump} gives
\[
d_{\bar g}\bigl(F(x_n,t_n),F(x_{n+1},t_{n+1})\bigr)
\leq R\sqrt{\tau_n}.
\]

Since $\sum_n\sqrt{\tau_n}<\infty$, the sequence $F(x_n,t_n)$
is Cauchy. By completeness of $M$, it converges to a point
$X_\infty\in M$, and
\[
d_{\bar g}\bigl(F(x_n,t_n),X_\infty\bigr)
\leq\sum_{m=n}^{\infty}R\sqrt{\tau_m}
=\frac{R}{1-K^{-1/2}}\sqrt{\tau_n}.
\]
By \eqref{eq:typeI-terminal},
$d_{\bar g}\bigl(F_T(x_n),F(x_n,t_n)\bigr)\to0$.
Thus $F_T(x_n)\to X_\infty$, and the compactness of
$F_T(\Sigma)$ gives $X_\infty\in F_T(\Sigma)$.
This proves \eqref{eq:typeI-essential-chain} with
$B=R/(1-K^{-1/2})$.
\end{proof}

\begin{proof}[Proof of Theorem~\ref{thm:no-type-I}]
Suppose that $|h|$ is unbounded on $\Sigma\times[0,T)$.
Then $T$ is the maximal time of smooth existence, since any smooth
extension across $T$ would give a uniform bound for $|h|$ on the
compact set $\Sigma\times[0,T]$. Under the Type~I bound
\eqref{eq:typeI-bound}, Lemma~\ref{lem:typeI-chain} gives constants
$\varepsilon,B>0$, a point $X_\infty\in F_T(\Sigma)$, and sequences
$t_n\uparrow T$ and $x_n\in\Sigma$ satisfying
\eqref{eq:typeI-essential-chain}.

The distance estimate in \eqref{eq:typeI-essential-chain} gives
$x_n\in N_{X_\infty}(t_n)$ for this choice of $B$. Applying
Lemma~\ref{lem:typeI-cone-decay} with the fixed center $X_\infty$,
we obtain
\[
\varepsilon
\leq (T-t_n)|h(x_n,t_n)|^2
\leq (T-t_n)\sup_{x\in N_{X_\infty}(t_n)}|h(x,t_n)|^2
\longrightarrow 0.
\]
This contradiction proves \eqref{eq:typeI-bounded}.
Theorem~\ref{thm:extension} now gives a smooth extension across $T$.
\end{proof}

This excludes finite-time Type~I singularities and proves
Theorem~\ref{thm:main}.

\begin{appendices}
\section{Curvature conventions and submanifold identities}
\label{app:submanifold-geometry}

We collect the curvature conventions and submanifold identities
used in Section~\ref{sec:evolution-formulae},
following Andrews and Baker~\cite{AB}.

We use the curvature convention
\[
R_\nabla(X,Y)\xi
:=\nabla_X\nabla_Y\xi-\nabla_Y\nabla_X\xi
  -\nabla_{[X,Y]}\xi.
\]
\subsection{Pullback bundles and connections}

Let $f:M\to N$ be a smooth map and $E$ a smooth vector bundle
over $N$. The pullback bundle $f^*E$ has fiber
$(f^*E)_x=E_{f(x)}$. For $\xi\in\Gamma(E)$, write
$\xi_f(x):=\xi(f(x))$ for the corresponding section of $f^*E$.
We recall the following standard properties.

\begin{proposition}\label{prop:curvature}
A connection $\nabla$ on $E$ induces a unique connection
${}^f\nabla$ on $f^*E$ satisfying
\[
({}^f\nabla_u\xi_f)_x
=(\nabla_{df_x(u)}\xi)_{f(x)}
\]
for $x\in M$, $u\in T_xM$, and $\xi\in\Gamma(E)$.
If $\nabla$ is compatible with a metric $g$ on $E$, then
${}^f\nabla$ is compatible with the pullback metric $g_f$,
defined by $(g_f)_x=g_{f(x)}$.
Its curvature satisfies
\[
R_{{}^f\nabla}(u,v)\eta
=R_\nabla(df_x(u),df_x(v))\eta
\]
for $u,v\in T_xM$ and $\eta\in(f^*E)_x$.
\end{proposition}

\begin{proposition}[Pullback torsion identity]
\label{prop:pullback-torsion}
If $\nabla$ is a torsion-free connection on $TN$, then
\[
{}^f\nabla_U(f_*V)-{}^f\nabla_V(f_*U)
=f_*[U,V]
\]
for all $U,V\in\Gamma(TM)$.
\end{proposition}

For an immersion $F:M^n\to N^{n+k}$ into a Riemannian manifold
$(N,\bar g)$, the differential induces an injective bundle map
\[
F_*:TM\to F^*TN,
\qquad
u\in T_xM\longmapsto(x,dF_x(u)).
\]
We identify $TM$ with its image under this map. Its orthogonal
complement with respect to the pullback metric $\bar g_F$ is
the normal bundle $N_FM$. The restrictions of $\bar g_F$ induce
metrics $g$ and $g^\perp$ on $TM$ and $N_FM$, respectively.

\subsection{Subbundles}

Let $E=K\oplus L$ be an orthogonal splitting of a metric vector
bundle $(E,g)$ over $M$, and let $\nabla$ be a metric-compatible
connection on $E$. Denote the inclusions by $\iota_K,\iota_L$
and the orthogonal projections by $\pi_K,\pi_L$.
The induced connection and second fundamental form of $K$ are
\[
\nabla^K_U\xi
:=\pi_K\bigl(\nabla_U(\iota_K\xi)\bigr),
\qquad
h^K(U,\xi)
:=\pi_L\bigl(\nabla_U(\iota_K\xi)\bigr),
\]
for $U\in\Gamma(TM)$ and $\xi\in\Gamma(K)$.
Define $\nabla^L$ and $h^L$ analogously.

Let $R^K$ denote the curvature of $\nabla^K$.
The Gauss and Codazzi equations are
\begin{equation}\label{eq:Gauss}
\begin{aligned}
R^K(U,V)\xi
={}&\pi_K\bigl(R_\nabla(U,V)(\iota_K\xi)\bigr)-h^L(U,h^K(V,\xi))+h^L(V,h^K(U,\xi)),
\end{aligned}
\end{equation}
and
\begin{equation}\label{eq:Codazzi-expanded}
\begin{aligned}
\pi_L\bigl(R_\nabla(U,V)(\iota_K\xi)\bigr)
={}&\nabla^L_U\bigl(h^K(V,\xi)\bigr)
-\nabla^L_V\bigl(h^K(U,\xi)\bigr)\\
&+h^K(U,\nabla^K_V\xi)-h^K(V,\nabla^K_U\xi)
-h^K([U,V],\xi),
\end{aligned}
\end{equation}
for $U,V\in\Gamma(TM)$ and $\xi\in\Gamma(K)$.
The induced connections preserve the restricted metrics
$g^K$ and $g^L$, and the Weingarten relation is
\begin{equation}\label{eq:weingarten}
g^L(h^K(U,\xi),\eta)+g^K(\xi,h^L(U,\eta))=0
\end{equation}
for $U\in\Gamma(TM)$, $\xi\in\Gamma(K)$, and
$\eta\in\Gamma(L)$.

For the immersion $F$ above, apply these constructions to
$F^*TN=TM\oplus N_FM$ with connection ${}^F\bar\nabla$,
where $\bar\nabla$ is the ambient Levi-Civita connection.
The induced connections on $TM$ and $N_FM$ are the
Levi-Civita connection $\nabla$ of $g$ and the normal
connection $\nabla^\perp$, respectively.
The tensor $h^{TM}$ is the second fundamental form of $F$.
The Gauss equation \eqref{eq:Gauss} applied to these two
subbundles gives the usual Gauss and Ricci equations.

\subsection{Time-dependent immersions}

Let $F:\Sigma^n\times I\to(N^{n+k},\bar g)$ be a smooth family
of immersions, and let $t:\Sigma\times I\to I$ be the time
projection. Set
\[
\mathcal H:=\ker dt,
\qquad
\mathcal N:=F_*(\mathcal H)^\perp\subset F^*TN,
\]
and identify $\mathcal H$ with $F_*(\mathcal H)$.
Denote the orthogonal projections by $\pi$ and $\pi^\perp$,
and the inclusion of $\mathcal N$ by $\iota$.
The pullback metric and connection induce metrics $g,g^\perp$
and metric connections $\nabla,\nabla^\perp$ on
$\mathcal H,\mathcal N$ in all spacetime directions.
Let $h^{\mathcal H}$ and $h^{\mathcal N}$ be the second
fundamental forms of these subbundles, and let $h$ be the
spatial restriction of $h^{\mathcal H}$.

Throughout, we assume that the variation is normal,
so $\pi F_*\partial_t=0$, as for \eqref{eq:gradient-flow}.
Proposition~\ref{prop:pullback-torsion} implies that $h$ is
symmetric. Since $[\partial_t,V]\in\Gamma(\mathcal H)$,
the same proposition gives
\begin{equation}\label{eq:h(i,t)}
h^{\mathcal H}(\partial_t,V)
=\pi^\perp\bigl({}^F\bar\nabla_V(F_*\partial_t)\bigr)
=\nabla^\perp_V(\pi^\perp F_*\partial_t)
\end{equation}
for $V\in\Gamma(\mathcal H)$.

For $U\in\Gamma(\mathcal H)$ and $\xi\in\Gamma(\mathcal N)$,
define the Weingarten map by
\[
\mathcal W(U,\xi)
:=-h^{\mathcal N}(U,\xi)
=-\pi\bigl({}^F\bar\nabla_U(\iota\xi)\bigr).
\]
By \eqref{eq:weingarten} and \eqref{eq:h(i,t)}, we have
\begin{align}
g^\perp(h(U,V),\xi)
&=g(V,\mathcal W(U,\xi)),
\label{eq:weingarten2}\\
g(h^{\mathcal N}(\partial_t,\xi),V)
&=-g^\perp\bigl(
\xi,\nabla^\perp_V(\pi^\perp F_*\partial_t)
\bigr),\notag
\end{align}
for $V\in\Gamma(\mathcal H)$.

Let $R$ and $R^\perp$ denote the curvatures of $\nabla$
and $\nabla^\perp$, respectively. We use the conventions
\[
\begin{aligned}
R(X,Y,W,Z)&=-g(R(X,Y)W,Z),\\
R^\perp(X,Y,\xi,\eta)
&=-g^\perp(R^\perp(X,Y)\xi,\eta),
\end{aligned}
\]
for $X,Y\in\Gamma(T(\Sigma\times I))$,
$W,Z\in\Gamma(\mathcal H)$, and
$\xi,\eta\in\Gamma(\mathcal N)$.
We use the same sign convention for the ambient curvature
tensor $\bar R$.

For $U,V,W,Z\in\Gamma(\mathcal H)$, applying
\eqref{eq:Gauss} to $\mathcal H$ gives
\begin{subequations}
\begin{align}
R(U,V)W
={}&-\mathcal W(V,h(U,W))+\mathcal W(U,h(V,W))
+\pi\bigl(\bar R(F_*U,F_*V)F_*W\bigr),
\label{eq:spatialGauss}\\
R(U,V,W,Z)
={}&g^\perp(h(U,W),h(V,Z))
-g^\perp(h(V,W),h(U,Z))
+(F^*\bar R)(U,V,W,Z).
\label{eq:spatialGauss2}
\end{align}
\end{subequations}
Applying \eqref{eq:Gauss} to $\mathcal N$ gives the Ricci equation
\begin{equation}\label{eq:Ricci}
R^\perp(U,V)\xi
=-h(V,\mathcal W(U,\xi))
+h(U,\mathcal W(V,\xi))
+\pi^\perp\bigl(\bar R(F_*U,F_*V)(\iota\xi)\bigr).
\end{equation}
The spatial Codazzi identity is
\begin{equation}\label{eq:spatialCodazzi}
(\nabla_Uh)(V,W)-(\nabla_Vh)(U,W)
=\pi^\perp\bigl(\bar R(F_*U,F_*V)F_*W\bigr).
\end{equation}

The connections $\nabla$ and $\nabla^\perp$ induce a connection
on $\mathcal H^*\otimes\mathcal H^*\otimes\mathcal N$.
In particular,
\[
(\nabla_{\partial_t}h)(V,W)
=\nabla^\perp_{\partial_t}\bigl(h(V,W)\bigr)
-h(\nabla_{\partial_t}V,W)-h(V,\nabla_{\partial_t}W).
\]
The corresponding mixed space--time Codazzi identity is
given by the following proposition.

\begin{proposition}\label{timelikeCodazzi}
Assume that the variation is normal. Then, for
\(V,W\in\Gamma(\mathcal H)\),
\begin{equation*}
\pi^\perp\bigl(\bar R(F_*\partial_t,F_*V)F_*W\bigr)
=(\nabla_{\partial_t}h)(V,W)
-(\nabla^{\perp 2}(\pi^\perp F_*\partial_t))(V,W)
-h\bigl(W,\mathcal W(V,\pi^\perp F_*\partial_t)\bigr),
\end{equation*}
where $(\nabla^{\perp 2}\xi)(V,W)
=\nabla^\perp_V\nabla^\perp_W\xi
-\nabla^\perp_{\nabla_VW}\xi$ for $\xi\in\Gamma(\mathcal N)$.
\end{proposition}

\begin{proof}
Taking $K=\mathcal H$ and $L=\mathcal N$ in
\eqref{eq:Codazzi-expanded} and using the definition of
$\nabla_{\partial_t}h$, we obtain
\begin{align*}
\pi^\perp\bigl(\bar R(F_*\partial_t,F_*V)F_*W\bigr)
={}&(\nabla_{\partial_t}h)(V,W)
-\nabla^\perp_V\bigl(h^{\mathcal H}(\partial_t,W)\bigr)
\notag\\
&+h^{\mathcal H}(\partial_t,\nabla_VW)
+h(\nabla_{\partial_t}V-[\partial_t,V],W).
\end{align*}
Since \(F\) is a normal variation, we have
\(\pi F_*\partial_t=0\). Thus
equation~\eqref{eq:h(i,t)} gives
\begin{align*}
&\nabla^\perp_V\bigl(h^{\mathcal H}(\partial_t,W)\bigr)
-h^{\mathcal H}(\partial_t,\nabla_VW)
\notag\\
&\quad=\nabla^\perp_V\nabla^\perp_W(\pi^\perp F_*\partial_t)
-\nabla^\perp_{\nabla_VW}(\pi^\perp F_*\partial_t)
\notag\\
&\quad=(\nabla^{\perp 2}(\pi^\perp F_*\partial_t))(V,W).
\end{align*}
The tangential projection of the pullback torsion identity gives
\[
\nabla_{\partial_t}V-[\partial_t,V]
=\pi\bigl({}^F\bar\nabla_V(F_*\partial_t)\bigr)
=-\mathcal W(V,\pi^\perp F_*\partial_t).
\]
Combining the preceding identities and using the symmetry of
\(h\) proves the result.
\end{proof}

For subsequent computations, we state without proof the following
identities for the covariant derivatives of the inclusion and
projection maps.

\begin{proposition}
With respect to the connections induced by
${}^F\bar\nabla$, $\nabla$, and $\nabla^\perp$, the bundle maps
$F_*\in\Gamma(\mathcal H^*\otimes F^*TN)$,
$\iota\in\Gamma(\mathcal N^*\otimes F^*TN)$,
$\pi\in\Gamma((F^*TN)^*\otimes\mathcal H)$, and
$\pi^\perp\in\Gamma((F^*TN)^*\otimes\mathcal N)$ satisfy
\begin{align}
(\nabla_U F_*)(V)
&=\iota h(U,V),\label{eq:nablaF*}\\
(\nabla_U\iota)(\xi)
&=-F_*\mathcal W(U,\xi),\notag\\
(\nabla_U\pi)(X)
&=\mathcal W(U,\pi^\perp X),\notag\\
(\nabla_U\pi^\perp)(X)
&=-h(U,\pi X),\label{eq:nablapiN}
\end{align}
for $U,V\in\Gamma(\mathcal H)$,
$\xi\in\Gamma(\mathcal N)$, and $X\in\Gamma(F^*TN)$.
\end{proposition}

We illustrate the preceding identities by deriving Simons' identity.
Here and below, we use the notation
$\nabla_W\nabla_Zh(U,V)
:=\bigl(\nabla_W(\nabla_Zh)-\nabla_{\nabla_WZ}h\bigr)(U,V)$
for the second covariant derivative of $h$.

\begin{proposition}[Simons' identity]\label{prop:Simons}
For $U,V,W,Z\in\Gamma(\mathcal H)$, we have
\begin{align*}
\nabla_W\nabla_Zh(U,V)-\nabla_U\nabla_Vh(W,Z)
={}&h(V,\mathcal W(U,h(W,Z)))
-h(Z,\mathcal W(W,h(U,V)))\\
&-h(U,\mathcal W(W,h(V,Z)))
+h(W,\mathcal W(U,h(V,Z)))\\
&+h(Z,\mathcal W(U,h(W,V)))
-h(V,\mathcal W(W,h(U,Z)))\\
&+\pi^\perp\bigl(\bar R(F_*W,F_*U)\iota h(V,Z)\bigr)\\
&-h(U,\pi(\bar R(F_*W,F_*V)F_*Z))
-h(W,\pi(\bar R(F_*Z,F_*U)F_*V))\\
&-h(Z,\pi(\bar R(F_*W,F_*U)F_*V))
-h(V,\pi(\bar R(F_*W,F_*U)F_*Z))\\
&+\pi^\perp\bigl(\bar R(\iota h(U,W),F_*V)F_*Z\bigr)
+\pi^\perp\bigl(\bar R(\iota h(W,Z),F_*U)F_*V\bigr)\\
&+\pi^\perp\bigl(\bar R(F_*W,\iota h(U,V))F_*Z\bigr)
+\pi^\perp\bigl(\bar R(F_*Z,\iota h(W,U))F_*V\bigr)\\
&+\pi^\perp\bigl(\bar R(F_*W,F_*V)\iota h(U,Z)\bigr)
+\pi^\perp\bigl(\bar R(F_*Z,F_*U)\iota h(W,V)\bigr)\\
&+\pi^\perp\bigl((\bar\nabla_{F_*U}\bar R)(F_*W,F_*V)F_*Z\bigr)
+\pi^\perp\bigl((\bar\nabla_{F_*W}\bar R)(F_*Z,F_*U)F_*V\bigr).
\end{align*}
\end{proposition}

\begin{proof}
Since the identity is tensorial, we choose local spatial
vector fields $U,V,W,Z$ whose spatial covariant derivatives
vanish at a fixed point. All computations below are
performed at that point. The Codazzi identity, the symmetry
of $h$, and the curvature commutation formula give
\begin{align*}
\nabla_W\nabla_Zh(U,V)
={}&\nabla^\perp_W\bigl((\nabla_Uh)(Z,V)\bigr)
+\nabla^\perp_W\bigl(\pi^\perp(\bar R(F_*Z,F_*U)F_*V)\bigr)\\
={}&\nabla_U\nabla_Wh(V,Z)+(R(W,U)h)(V,Z)\\
&+\nabla^\perp_W\bigl(\pi^\perp(\bar R(F_*Z,F_*U)F_*V)\bigr)\\
={}&\nabla_U\nabla_Vh(W,Z)+(R(W,U)h)(V,Z)\\
&+\nabla^\perp_W\bigl(\pi^\perp(\bar R(F_*Z,F_*U)F_*V)\bigr)
+\nabla^\perp_U\bigl(\pi^\perp(\bar R(F_*W,F_*V)F_*Z)\bigr).
\end{align*}
Here the Codazzi identity is used in the first and third
equalities, and the curvature commutation formula in the
second. Since $h$ is an $\mathcal N$-valued bilinear form
on $\mathcal H$, the induced curvature acts by
\[
(R(W,U)h)(V,Z)
=R^\perp(W,U)h(V,Z)
-h(R(W,U)V,Z)-h(V,R(W,U)Z).
\]
This expression is expanded using the Gauss equation
\eqref{eq:spatialGauss} and the Ricci equation \eqref{eq:Ricci}.

For the first ambient curvature derivative term, the
product rule gives
\begin{align*}
\nabla^\perp_W\bigl(\pi^\perp(\bar R(F_*Z,F_*U)F_*V)\bigr)
={}&(\nabla_W\pi^\perp)\bigl(\bar R(F_*Z,F_*U)F_*V\bigr)
+\pi^\perp\bigl({}^F\bar\nabla_W
  (\bar R(F_*Z,F_*U)F_*V)\bigr).
\end{align*}
By \eqref{eq:nablapiN}, the first term is
$-h(W,\pi(\bar R(F_*Z,F_*U)F_*V))$.
For the second term, the definition of the pullback
connection and \eqref{eq:nablaF*} yield
\begin{align*}
{}^F\bar\nabla_W(\bar R(F_*Z,F_*U)F_*V)
={}&(\bar\nabla_{F_*W}\bar R)(F_*Z,F_*U)F_*V\\
&+\bar R(\iota h(W,Z),F_*U)F_*V
+\bar R(F_*Z,\iota h(W,U))F_*V\\
&+\bar R(F_*Z,F_*U)\iota h(W,V).
\end{align*}
The remaining term
$\nabla^\perp_U\bigl(\pi^\perp(\bar R(F_*W,F_*V)F_*Z)\bigr)$
is expanded in the same way. Substituting these expansions
and using the symmetry of $h$ proves the identity.
\end{proof}

\begin{remark}[Traced Simons identity]
\label{rem:traced-Simons}
Let $\{e_k\}$ be a local orthonormal frame of
$\mathcal H$. Taking $W=Z=e_k$ in
Proposition~\ref{prop:Simons} and summing over $k$ yields
\begin{equation}\label{traced-Simons-identity}
\Delta h
=\nabla^{\perp2}\HH
+\mathcal A_{\HH}
+\mathcal R_{\HH}
+\mathcal S(h,\bar R,\bar\nabla\bar R),
\end{equation}
where
\[
\mathcal A_{\HH}(U,V)
=h\bigl(V,\mathcal W(U,\HH)\bigr),
\qquad
\mathcal R_{\HH}(U,V)
=\pi^\perp\bigl(
\bar R(\iota\HH,F_*U)F_*V
\bigr).
\]
Here \(\mathcal S(h,\bar R,\bar\nabla\bar R)\) denotes
the sum of all remaining terms in the trace. In particular,
it contains no covariant derivatives of \(h\).
\end{remark}
\end{appendices}

\section{Acknowledgments} I am grateful to my advisor, Professor Jiayu Li, for his guidance and support throughout this work.

\end{document}